\documentclass[a4paper,11pt,reqno]{amsart}
\usepackage[centertags]{amsmath}
\usepackage{amsfonts,amssymb,amsthm,dsfont,cases,amscd,esint,enumerate,fancyhdr,newlfont,color}
\usepackage[english]{babel}
\usepackage[body={15cm,21.5cm},centering]{geometry} 
\newtheorem{theor}{Theorem}[section]
\newtheorem{lem}[theor]{Lemma}

\newtheorem{cor}[theor]{Corollary}
\theoremstyle{definition}

\newtheorem{rem}[theor]{Remark}

\numberwithin{equation}{section}

\newcommand{\N}{\mathbb N}

\newcommand{\Z}{\mathbb Z}
\newcommand{\e}{\varepsilon}
\newcommand{\dist}{\operatorname{dist}}
\newcommand{\R}{\mathbb R}

\newcommand{\Kc}{\mathcal K}
\newcommand{\Sc}{\mathcal S}

\newcommand{\loc}{{\operatorname{loc}}}
\newcommand{\Id}{\operatorname{Id}}
\newcommand{\E}{\mathbb{E}}

\newcommand{\Ld}{\operatorname{L}}

\newcommand{\Pm}{\mathbb{P}}
\newcommand{\expec}[1]{\mathbb{E}\left[ #1 \right]}
\newcommand{\expecm}[1]{\mathbb{E}\big[ #1 \big]}

\newcommand{\step}[1]{\noindent \textit{Step} #1.}
\newcommand{\substep}[1]{\noindent \textit{Substep} #1.}

\usepackage{hyperref}

\title[Lower spectrum of heterogeneous acoustic operators]{Large-scale dispersive estimates for acoustic operators: homogenization meets localization}

\author[M. Duerinckx]{Mitia Duerinckx}
\address[Mitia Duerinckx]{Universit\'e Libre de Bruxelles, D\'epartement de Math\'ematique, 1050~Brussels, Belgium}
\email{mitia.duerinckx@ulb.be}
\author[A. Gloria]{Antoine Gloria}
\address[Antoine Gloria]{Sorbonne Universit\'e, CNRS, Universit\'e de Paris, Laboratoire Jacques-Louis Lions, 75005~Paris, France \& Institut Universitaire de France \& Universit\'e Libre de Bruxelles, D\'epartement de Math\'ematique, 1050~Brussels, Belgium}
\email{antoine.gloria@sorbonne-universite.fr}

\begin{document}

\begin{abstract}
This work relates quantitatively homogenization to Anderson localization for acoustic operators in disordered media.
By blending dispersive estimates for homogenized operators and quantitative homogenization of the wave equation, we derive large-scale dispersive estimates for waves in disordered media that we apply to the spreading of low-energy eigenstates.
This gives a short and direct proof 
that the lower spectrum of the acoustic operator is purely absolutely continuous in case of periodic
media,
and it further provides new lower bounds on the localization length of possible eigenstates in case of quasiperiodic or random
media.
\end{abstract}


\maketitle
\setcounter{tocdepth}{1}
\tableofcontents

\section{Introduction and main results}

\subsection{General overview}
Given a coefficient field $A:\R^d\to\R^{d\times d}$ that is symmetric and uniformly elliptic in the sense of
\begin{equation}\label{eq:ell-cond}
e\cdot A(x)e\ge\tfrac1{C_0}|e|^2,\qquad|A(x)e|\le |e|,\qquad\text{for all $x,e\in\R^d$},
\end{equation}
for some $C_0\ge 1$, we consider the associated acoustic operator
\[H_a:=-\nabla\cdot A\nabla,\qquad \text{on $\Ld^2(\R^d)$}.\]
 Depending on properties of $A$, the spectrum of $H_a$ can be of different types, but not much is known in general: there has indeed been surprisingly few studies on the spectrum of this operator despite its role in the description of classical wave propagation in disordered media. Only two particular settings are well understood:%
\begin{enumerate}[---]
\item If the coefficient field $A$ is periodic, then the lower spectrum of $H_a$ is absolutely continuous in any dimension $d\ge1$ (see e.g.~\cite[Section~6.3]{Kuchment-16}).
In fact, the whole spectrum of $H_a$ is conjectured to be continuous in that case provided that~$A$ is smooth enough (see~\cite[Conjecture~6.13]{Kuchment-16}), while large eigenvalues are known to exist in general in the non-smooth setting~\cite{Filonov}.
\smallskip\item If the coefficient field $A$ is random (at least, in case of the so-called random displacement model defined in~\cite[eqn~(1.3)]{Sims-Stolz}), and if the space dimension is~$d=1$, then the operator~$H_a$ {displays Anderson localization in the sense that it} has dense pure point spectrum with exponentially decaying eigenstates~\cite{Sims-Stolz}  --- just like $1$-dimensional random Schr\"odinger operators.
\end{enumerate}
In the higher-dimensional random setting,
the comparison with random Schr\"odinger operators is less clear: for $d>2$, while Schr\"odinger operators are expected to display a metal--insulator transition, there is still no definite consensus on what to expect for the type of the lower spectrum of $H_a$.
The results in~\cite{Figotin-Klein} show that Anderson localization holds for $H_a$ near all spectral edges except near the bottom of the spectrum, where the situation remains unclear.
The particularity of the lower spectrum of $H_a$ originates in its degeneracy:
the constant function~$1$ is always an extended state at energy~$0=\inf\sigma(H_a)$, in link with the expected degeneracy of Lyapunov exponents, cf.~\cite{Sims-Stolz}.

In the present contribution, we focus on the spectrum of~$H_a$ close to the lowest energy~$0$ and we show that specific techniques such as homogenization theory 
can be pushed forward to bring valuable new information in that spectral region, cf.~Theorem~\ref{th:main} below.
This is nontrivial as homogenization theory is a priori far from giving any information on the resolvent of~$H_a$ close to its spectrum. Instead of elliptic homogenization, the crux of our approach is to {build upon} recent results on the long-time homogenization for the wave equation associated with~$H_a$. More precisely, these hyperbolic homogenization results are used to first derive {what we call} large-scale dispersive estimates, cf.~Theorem~\ref{th:disp} below, from which we deduce strong estimates on the spatial spreading of eigenstates if any.

\subsection{Main result: spreading of eigenstates}
Our main result below provides a strong control on  the spatial spreading of the mass density of eigenstates in the lower spectrum (as well as on the spreading of their energy density, cf.~Remark~\ref{rem:supp}). This is expressed as a lower bound on the spatial `width', or localization length, of an eigenstate~$\psi$ if it exists, which we define as
\begin{equation}\label{eq:width}
\ell_\theta(\psi)\,:=\,\inf\big\{r\ge0:\|\psi\|_{\Ld^2(B_r)}\ge{(1-\theta)\|\psi\|_{\Ld^2(\R^d)}}\big\},\qquad\text{given $0<\theta<\tfrac12$}.
\end{equation}
Our approach is quite robust and applies to essentially any kind of disorder for the coefficient field $A$ under suitable statistical spatial homogeneity assumptions, which we illustrate by considering at the same time periodic, quasiperiodic, and random coefficient fields in any dimension. In the random setting, we only cover the strongest mixing conditions for {conciseness}.\footnote{Under weaker mixing conditions, we would be led to weaker lower bounds on the localization length, and we refer to~\cite{GNO-quant,DGR} for the corresponding homogenization tools.}
 We emphasize that the spreading of low-energy eigenstates is particularly striking when compared to generic low-energy functions obtained by rescaling: indeed, 
for~$\phi_\lambda:=\lambda^{d/4}\phi(\lambda^{1/2}\cdot)$,
we find
\begin{equation}\label{e.discus-sc}
\|(H_a)^\frac12\phi_\lambda\|_{\Ld^2(\R^d)}\simeq
\lambda^{\frac12}\|\phi\|_{\Ld^2(\R^d)},\qquad\ell_\theta(\phi_\lambda)=\lambda^{-\frac12}\ell_\theta(\phi),
\end{equation}
while, in the result below, any possible eigenstate $\psi_\lambda$ at energy $\lambda$ is shown to have a width~$\ell_\theta(\psi_\lambda)$ always substantially much larger than $O(\lambda^{-1/2})$.
The proof is given in Section~\ref{sec:th} and is based on long-time hyperbolic homogenization in form of the large-scale dispersion estimates that {we introduce} in Section~\ref{sec:dispersion-intro} below.

\begin{theor}\label{th:main}$ $
Let $d\ge1$ and let the coefficient field $A:\R^d\to\R^{d\times d}$ be symmetric and uniformly elliptic in the sense of~\eqref{eq:ell-cond}.
\begin{enumerate}[(i)]
\item \emph{Periodic setting:}\\
Assume that $A$ is periodic. Then there exists $\lambda_0>0$ (only depending on $d,C_0$) such that the operator $H_a$ admits no eigenvalue in $[0,\lambda_0]$.
\smallskip\item \emph{Quasiperiodic setting:}\\
Assume that $A$ takes the form
$A(x)=A_0(Fx)$
for some frequency matrix $F\in\R^{M\times d}$, with $M>d$, and for some lifted map $A_0:\R^M\to\R^{d\times d}$ that is periodic on~$\R^M$ and Gevrey regular. Further assume that the frequency matrix $F$ satisfies a Diophantine condition, that is, for some $C,\gamma>0$,
\[\qquad|Fz|\,\ge\,\tfrac1C|z|^{-\gamma},\qquad\text{for all $z\in\Z^M\setminus\{0\}$}.\]
Then, given $0<\theta<\frac12$, there exists {$\lambda_{\theta}>0$} with the following property:
if~$H_a$ has an eigenvalue $\lambda\le \lambda_\theta$, any associated eigenstate~$\psi_\lambda$ satisfies
\begin{equation*}
\qquad\ell_\theta(\psi_\lambda)\,\ge\,\exp(\lambda^{-\sigma}),
\end{equation*}
 where we recall the notation~\eqref{eq:width} for the spatial spreading of $\psi_\lambda$, and where~$\sigma>0$ only depends on the Diophantine exponent~$\gamma$ and on the Gevrey regularity of $A_0$.
\smallskip\item \emph{Random setting:}\\
 Assume that $A=\tilde A(\cdot,\omega_0)$ is a realization of a random field $\tilde A$: more precisely, consider a probability space $(\Omega,\Pm)$ and a statistically spatially homogeneous, uniformly elliptic, random field $\tilde A:\Omega\times\R^d\to\R^{d\times d}$, in the following sense,
\begin{enumerate}[---]
\smallskip\item \emph{Random field:} $\tilde A$ is a jointly measurable map $\R^d\times\Omega\to\R^{d\times d}$ such that for all $x\in\R^d$ the function $\tilde A(x,\cdot):\Omega\to\R^{d\times d}$ is measurable.
\smallskip\item \emph{Symmetry and uniform ellipticity:} for $\Pm$-almost all $\omega$, the realization $\tilde A(\cdot,\omega)$ is a symmetric matrix field and satisfies the uniform ellipticity condition~\eqref{eq:ell-cond}.
\smallskip\item \emph{Statistical spatial homogeneity:} the finite-dimensional law of $\tilde A$ is shift-invariant, that is, for any finite set $E\subset\R^d$ the law of $\{\tilde A(x+y,\cdot)\}_{x\in E}$ does not depend on the shift $y\in\R^d$.
\end{enumerate}
\smallskip
Further assume that $\tilde A$ satisfies strong enough mixing assumptions, for instance:\footnote{ The same result holds if, instead of the finite range of dependence, we rather assume for instance that $\tilde A$ has a Gaussian structure with integrable covariance --- more precisely, this means to assume that~$\tilde A$ can be written as $\tilde A(x,\omega)=A_0(G(x,\omega))$, where $G:\R^d\times\Omega\to\R^p$ is a statistically spatially homogeneous Gaussian field, for some $p\ge1$, where $A_0\in C^1_b(\R^p;\R^{d\times d})$, and where we require the covariance function of~$G$ to be written as $\expec{G(x,\cdot)\otimes G(y,\cdot)}=(C_0\ast C_0)(x-y)$ for some even convolution kernel $C_0:\R^d\to\R^{p\times p}$ satisfying the integrability condition $\int_{\R^d}(\sup_{B(x)}|C_0|)\,dx<\infty$.}
\begin{enumerate}[---]
\smallskip\item \emph{Finite range of dependence:} for all finite subsets $E,F\subset\R^d$ with $\dist(E,F)\ge1$, the families of random variables $(\tilde A(x,\cdot))_{x\in E}$ and $(\tilde A(x,\cdot))_{x\in F}$ are independent.
\end{enumerate}
\smallskip
Then, given $0<\theta<\frac12$ and $\e>0$, there exists a positive random variable $\lambda_{\theta,\e}$ with finite (inverse) stretched exponential moments,
\[\qquad\E\Big[{\exp\big(\tfrac1C(\lambda_{\theta,\e})^{-\frac1C}\big)}\Big]<\infty,\qquad\text{for some $C<\infty$ (depending on $d,C_0,\theta,\e$)},\]
such that the following property holds:  for $\Pm$-almost all $\omega_0$, if for the coefficient field~$A=\tilde A(\cdot,\omega_0)$ the operator~$H_a$ has an eigenvalue $\lambda\le\lambda_{\theta,\e}(\omega_0)$, any associated eigenstate~$\psi_\lambda$ satisfies 
\[
\ell_\theta(\psi_\lambda)\,\ge\, 
\lambda^{\e-\frac12(\frac d2+1)}.
\]
\end{enumerate}
\end{theor}

\begin{rem}[Energy distribution]\label{rem:supp}
In the above result, in the definition~\eqref{eq:width} of the width of eigenstates, the local mass $\|\psi\|_{\Ld^2(B_r)}$ can be replaced by the local (rescaled) energy: more precisely, for an eigenstate $\psi_\lambda$ at energy $\lambda$, the same estimates hold as above if~$\ell_\theta(\psi_\lambda)$ is replaced by
\begin{equation}\label{eq:width-bis}
 \ell_\theta'(\psi_\lambda)\,:=\,\inf\Big\{r\ge0:\tfrac1{\sqrt{\lambda}} \|\sqrt A \nabla \psi_\lambda\|_{\Ld^2(B_r)} \ge(1-\theta)\|\psi_\lambda\|_{\Ld^2(\R^d)}\Big\}.
\end{equation}
Note that the rescaling of the local energy  is to ensure
\[\tfrac1{\sqrt{\lambda}} \|\sqrt A \nabla \psi_\lambda\|_{\Ld^2(B_r)}\,\xrightarrow{r\uparrow\infty}\,\big(\tfrac1{{\lambda}}\textstyle\int_{\R^d} \nabla \psi_\lambda \cdot A \nabla \psi_\lambda\big)^\frac12\,=\,\|\psi_\lambda\|_{\Ld^2(\R^d)}.\]
In other words, we show that the spatial spreading of the mass density and of the energy density of eigenstates are comparable. A proof is included in Step~7 of Section~\ref{sec:th}.
\end{rem}
We briefly comment on the above result and compare it to the literature, separately discussing the periodic, quasi-periodic, and random settings:

\begin{enumerate}[---]
\item {\it Periodic setting:}\\
In the periodic setting, item~(i) is well known to the expert and is part of the folklore of Floquet theory, see e.g.~\cite[Section~6.3]{Kuchment-16}.
In~\cite{AKS-20}, Armstrong, Kuusi, and Smart gave a different proof of~(i), based on periodic elliptic homogenization theory in form of a large-scale analyticity result for $A$-harmonic functions,
which prevents the existence of \emph{exponentially-localized} eigenstates in the lower spectrum.
To extend this to the non-existence of general~$\Ld^2$ eigenstates, the authors relied on a deep result by Kuchment based on Floquet theory, which states that the existence of an $\Ld^2$ eigenstate would imply the existence of an exponentially-localized eigenstate  (see~\cite[Theorem~6.15]{Kuchment-16}). 
In the present work, our proof of item~(i) also relies on homogenization theory, but now applied to the wave equation rather than to the elliptic equation associated with the operator~$H_a$. As opposed to~\cite{AKS-20}, the non-existence of general $\Ld^2$ eigenstates follows directly from our arguments without relying on Floquet theory.\footnote{Our proof would in fact be substantially simpler if we restricted ourselves to the study of exponentially-localized eigenstates.}
This latter point is crucial as it precisely allows us to extend our argument to more general disordered media, such as the quasiperiodic and random settings, which could not be covered so far by other available techniques.
{In the same survey article~\cite{Kuchment-16}, Kuchment conjectures that there is no $\Ld^2$ eigenstate even higher in the spectrum provided $A$ is periodic and smooth enough --- such a result is beyond the reach of arguments based on homogenization (which, by definition, focuses on low frequencies).}

\smallskip
\item {\it Quasiperiodic and random settings:}\\
The results of items~(ii) and~(iii) in the quasiperiodic and random settings are  new and provide lower bounds on the  spatial spreading of the mass density of potential eigenstates in the lower spectrum. This is new even in the~1D random setting of~\cite{Sims-Stolz}, where only the qualitative result $\ell_\theta(\psi_\lambda)\uparrow\infty$ was known so far as~$\lambda\downarrow0$.
As stated above, our result further provides a lower bound
which reads {$\ell_\theta(\psi_\lambda)\ge \lambda^{\e-3/4}$} in the 1D random setting as~$\lambda\downarrow0$ for any $0<\theta<\frac12$ and $\e>0$ (see also Corollary~\ref{cor-hyp} for an improved estimate in the specific setting of~\cite{Sims-Stolz}).

\smallskip\noindent
Yet, we do not expect these convergence rates to be sharp. Indeed, in the quasiperiodic setting, we do not expect eigenvalues to exist at all in the lower spectrum: the strong spatial spreading stated in item~(ii) should be seen as a partial indication in that direction --- while the continuity of the lower spectrum remains an open problem. Note that in the quasiperiodic setting the story is still incomplete even for Schr\"odinger operators beyond the 1D case of~\cite{Eliasson-92}:
the existence of continuous spectrum has been established in~\cite{Karpeshina-Shterenberg-14} in dimension~$d=2$ only for very specific quasiperiodic potentials, without even ruling out the existence of localized eigenstates. {In \cite{DGS-21}, Shirley and the authors construct formal series for extended states, which somehow correspond to the corrector expansions that we consider in the sequel, but these are not absolutely convergent; the exact same difficulty arises for the long-time homogenization of acoustic waves in~\cite{DGR}.}

\smallskip\noindent
As already mentioned, the problem is richer and even more poorly understood in the random setting, to the point that no definite conjecture {has been} formulated for the type of the lower spectrum of the operator~$H_a$ beyond the 1D case of~\cite{Sims-Stolz}. Item~(iii) {does provide} some new, nontrivial information, although {it is not expected to be} optimal.

\smallskip
\item {\it Some improvements in the random setting:}\\
The result of item~(iii)
holds for any random coefficient field that satisfies some form of independence,
and in general it is the best result that can be obtained from our method.
Yet, this does not mean that additional assumptions cannot yield an improved scaling, and the 1D random displacement model of~\cite{Sims-Stolz} turns out to be such an example. Indeed, the coefficient field $A$ in that case is such that $A^{-1}$ is hyperuniform in the sense of~\cite{Torquato-Stillinger-03}, which means that the variance of large-scale spatial averages of $A^{-1}$ decays faster than the CLT scaling.
In 1D, this is obviously transferred to averages of the first corrector gradient in homogenization theory (which is indeed proportional to~$A^{-1}-\E[A^{-1}]$),
thus implying that the corrector for this specific hyperuniform 1D model has the same growth as a generic 2D corrector; see discussion in~\cite[Section~2.7]{MR4366081}.
Hence, the proof of Theorem~\ref{th:main} yields the following improvement of Theorem~\ref{th:main}.

\begin{cor}\label{cor-hyp}
Let $\tilde A$ be given by the 1D random displacement model of~\cite[eqn~(1.3)]{Sims-Stolz}.
Given {$0<\theta<\frac12$} and $\e>0$, there exists a positive random variable $\lambda_{\theta,\e}$ with finite (inverse) stretched exponential moments with the following property:  for $\Pm$-almost all~$\omega_0$, if for the coefficient field~$A=\tilde A(\cdot,\omega_0)$ the operator~$H_a$ has an eigenvalue $\lambda\le\lambda_{\theta,\e}(\omega_0)$, any associated eigenstate~$\psi_\lambda$ satisfies 
\[\ell_\theta(\psi_\lambda)\,\ge\, \lambda^{\e-1}.\]
\end{cor}

\noindent
Hyperuniformity is not the only way to improve scalings: this issue is well-known to the expert in fluctuations in stochastic homogenization and is related in 1D to the Hermite rank, see~\cite{Taqqu,BGMP-08,LNZH-17,MR4414703}.
In any case, this observation suggests
that the localization length might depend on further specific features of the coefficient field, and that it could be for instance significantly larger for the random displacement model than for amorphous models such as Poisson models in 1D.
\end{enumerate}

\subsection{General strategy: homogenization meets localization}
We briefly describe the general idea
to use homogenization to obtain information on localized eigenstates. 
As opposed to previous works based on spectral-theoretical methods, or on the study of the elliptic eigenvalue problem $H_a \psi_\lambda=\lambda \psi_\lambda$, we adopt a dynamical point of view and capitalize on recent advances on the long-time homogenization of heterogeneous wave equations~\cite{BLP-78,BG,ALR,DGR}.  
{Let us first present the heuristic of our strategy (which, unfortunately, we cannot implement so straightforwardly).}
Our argument is the following: if~$H_a$ admits an eigenvalue, we can consider the associated standing wave, which consists of time oscillations around the eigenstate and solves the heterogeneous wave equation. If the eigenvalue is close to~$0$, homogenization theory can be applied and ensures that the low-energy standing wave must behave on some long time window like the solution of a corresponding `homogenized' constant-coefficient wave equation: in particular, it must display ballistic transport as the latter on this long time window.
{This should imply} that the eigenstate cannot be too localized: indeed, if the mass
of the eigenstate were supported in a ball that was too small, then the ballistic transport inherited by homogenization would end up pushing the standing wave away from its support --- a contradiction.
With this in mind, the statement of Theorem~\ref{th:main} is not surprising and can be expected to follow from long-time homogenization for waves.
For periodic coefficients, as homogenization holds up to arbitrarily long times at low energy, the argument {should imply that no eigenvalue can exist} in the lower spectrum, thus proving item~(i) in Theorem~\ref{th:main}. For quasiperiodic and random coefficients, homogenization ultimately breaks down on very long times: this leaves the possibility of existence of eigenstates in the lower spectrum, but we {should still be able to} deduce nontrivial lower bounds on the spatial spreading of low-energy eigenstates as stated in items~(ii) and~(iii).

Although fairly natural, the above general strategy of using long-time homogenization in form of transport properties {fails to prove} Theorem~\ref{th:main}.
{Indeed, this approach requires moment bounds on the eigenstate, which we do not assume and which cannot be fixed by an approximation argument.}
Instead, we shall rely on another important property of homogenized wave equations, that is, dispersion.
More precisely, rather than describing how far the wave gets transported, we shall use dispersion to show that the wave vanishes away in the supremum norm: if the mass of the standing wave ends up getting too small on the support of the considered eigenstate, we still get a contradiction.
The question is thus how to use homogenization to show that low-energy waves in disordered media indeed have similar dispersive properties as waves in homogeneous media on some long time window.
 This leads us to establishing {what we call} ``large-scale dispersive estimates'' as described in Section~\ref{sec:dispersion-intro} below,
which should be put on the same footing as the large-scale elliptic regularity theory recently developed for heterogeneous elliptic operators~\cite{AS,AKM-book,MR4103433} (see also the above-mentioned large-scale analyticity in~\cite{AKS-20}).

Another difficulty
is that potential eigenstates taken as initial data for the wave equation in our argument need not have a definite scaling in general:
they are oscillatory and do certainly not take the simple form of smooth rescaled bumps as considered in previous long-time homogenization results~\cite{BG,ALR,DGR} (such a form would indeed yield suboptimal results, cf.~\eqref{e.discus-sc}). Their oscillations thus need to be further monitored and controlled, preventing us from blindly appealing to the long-time homogenization results of~\cite{BG,ALR,DGR}, see Section~\ref{sec:corr}. This is taken care of by a suitable approximation argument, see Step~1 in Section~\ref{sec:disp}.

As lower bounds on the {localization length} of eigenstates are deduced from the competition between {dispersion} and localization,
we are led to a natural comment regarding the expected suboptimality of the results obtained with this general strategy:
beyond the time window for ballistic transport {and dispersion}, one actually expects the existence of another time window with only diffusive transport. Improved lower bounds on the spatial width of eigenstates might then follow from the competition between this diffusive transport and localization, {provided one can derive dispersive estimates from diffusive transport}. However, diffusive transport has not been established yet in any form for low-energy waves {(neither its relation to dispersive estimates)}, and remains an open problem.

\medskip

\subsection*{Notation}
\begin{enumerate}[---]
\item We denote by $C\ge1$ any constant that only depends on the dimension $d$ and on the ellipticity constant $C_0$ in~\eqref{eq:ell-cond}. We use the notation $\lesssim$ (resp.~$\gtrsim$) for $\le C\times$ (resp.~$\ge\frac1C\times$) up to such a multiplicative constant $C$. We write $\ll$ (resp.~$\gg$) for $\le C\times$ (resp.~$\ge C\times$) up to a sufficiently large multiplicative constant $C$. We add subscripts to indicate dependence of constants on other parameters.
\item The ball centered at $x$ of radius $r$ in $\R^d$ is denoted by $B_r(x)$, and we set $B(x)=B_1(x)$, $B_r=B_r(0)$, and $B=B_1(0)$.
\item Recall that $ \lfloor r \rfloor$ denotes the integer part of $r$, that is, the largest integer smaller or equal to $r$, and we also denote by $ \lceil r \rceil$ the smallest integer larger or equal to $r$.
\item We set $\langle s\rangle:=(1+|s|^2)^{1/2}$.
\item $\N$ stands for the set of positive integers.
\item {Given $v\in C^\infty_c(\R^d)$, we denote by $\hat v(\xi)=\int_{\R^d} e^{-ix\cdot \xi} v(x) dx$ its Fourier transform, 
and write the inversion formula as $v(x)=\int_{\R^d} e^{ix\cdot \xi} \hat v(\xi)d^*\xi$ with $d^*\xi=(2\pi)^{-d} d\xi$.}
\end{enumerate}

\subsection{Large-scale dispersive estimates}\label{sec:dispersion-intro}
We start by recalling the standard form of dispersive, i.e. $\Ld^1(\R^d)\to\Ld^\infty(\R^d)$, estimates for classical waves in an homogeneous medium: in any dimension $d\ge1$, given an initial condition $u^\circ\in\Ld^1\cap H^1(\R^d)$ that has Fourier transform supported in the unit ball $B$ (say), we have for all $t\ge0$,
\begin{equation}\label{eq:standard-disp-1}
\|\!\cos(t\sqrt{-\triangle})u^\circ\|_{\Ld^\infty(\R^d)}\,\lesssim\,(1+t)^{-\frac{d-1}2}\|u^\circ\|_{\Ld^1(\R^d)}.
\end{equation}
In case of localized data, a better time decay can actually be obtained for the supremum norm on the localized initial support when time exceeds the size of the support: more precisely, for any $\eta>0$, we further have for all $R,L\ge1$ and $t\gg R+L$,
\begin{equation}\label{eq:standard-disp-2}
\|\!\cos(t\sqrt{-\triangle})(u^\circ\mathds1_{B_L})\|_{\Ld^\infty(B_R)}\,\lesssim_\eta\,(1+t)^{\eta-d}\|u^\circ\|_{\Ld^1(B_L)}.
\end{equation}
This second estimate easily follows from the standard proof of~\eqref{eq:standard-disp-1}, based on a careful analysis of oscillatory Fourier integrals, see e.g.~\cite[Section~8.1.3]{BCD-11}.

In the spirit of the large-scale elliptic regularity inherited by heterogeneous elliptic equations by homogenization~\cite{AKM-book,MR4103433},
we may naturally {investigate} what kind of dispersion properties {can be proved for} waves in heterogeneous media.
As localized eigenvalues may appear in general for heterogeneous operators even in the lower spectrum, such as for $H_a$ in the 1D random setting~\cite{Sims-Stolz}, dispersive estimates~\eqref{eq:standard-disp-1} and~\eqref{eq:standard-disp-2} {cannot hold} in the same form: eigenvalues lead to time-periodic solutions that do not satisfy any time decay.

To {motivate} the form of {the large-scale dispersive estimates that we shall state} below,
we first comment on the periodic setting.
If $A$ is a periodic uniformly elliptic coefficient field, we may expect to prove dispersion by analyzing the integral representation of the wave by Floquet--Bloch instead of Fourier theory. Yet, the operator $H_a$ on $\Ld^2(\R^d)$ may then have eigenvalues in the upper spectrum in general, cf.~\cite{Filonov}, which would destroy dispersion: a projection $P_c$ on the continuous part of the spectrum is necessary and we only expect for all~$t\ge0$, provided that $u^\circ$ has again Fourier transform supported in~$B$ (say),
\begin{equation}\label{eq:periodic-disp-Pc}
\|\!\cos(t\sqrt{H_a})P_cu^\circ\|_{\Ld^\infty(\R^d)}\,\lesssim\,(1+t)^{-\frac{d-1}2}\|u^\circ\|_{\Ld^1(\R^d)}.
\end{equation}
Such an estimate has two main drawbacks in our context of general disordered media: it is an empty statement if the lower spectrum is pure point (which is already the case in the 1D random setting~\cite{Sims-Stolz}), and moreover the projection~$P_c$ {is not explicit and} depends on the disorder.
We thus need to further reformulate this estimate.
In the periodic setting, as at least the lower spectrum of $H_a$ is known to be purely absolutely continuous (see e.g.~\cite[Section~6.3]{Kuchment-16}, or also~Theorem~\ref{th:main}(i) above), the projection can be neglected in~$\Ld^2(\R^d)$ up to an error estimated by the energy norm: if the spectrum is continuous below $\lambda_0>0$, we find
\[\|(1-P_c)u\|_{\Ld^2(\R^d)}\,\le\,\lambda_0^{-\frac12}\|\sqrt{H_a} u\|_{\Ld^2(\R^d)}\,\lesssim\,\lambda_0^{-\frac12}\|\nabla u\|_{\Ld^2(\R^d)}.\]
 Averaging the above dispersive estimate~\eqref{eq:periodic-disp-Pc} on a ball of size $R>0$, and noting that the assumption that the Fourier transform of $u^\circ$ is supported in $B$ can also be removed up to an additional error term involving the energy norm,
we are led to the following estimate without projection: for all $t\ge0$, $R>0$, and $u^\circ\in\Ld^1\cap H^1(\R^d)$,
\begin{multline}\label{eq:periodic-disp-0}
\sup_{x\in\R^d}{R^{-\frac d2}}\|\!\cos(t\sqrt{H_a})u^\circ\|_{\Ld^2({B_R}(x))}
\,\lesssim\,(1+t)^{-\frac{d-1}2}\|u^\circ\|_{\Ld^1(\R^d)}\\[-1mm]
+(1\wedge\lambda_0)^{-\frac12}R^{-\frac d2}\|\nabla u^\circ\|_{\Ld^2(\R^d)}.
\end{multline}
Whereas \eqref{eq:periodic-disp-Pc} is a neat $\Ld^\infty-\Ld^1$ estimate, this estimate~\eqref{eq:periodic-disp-0} takes the form of a large-scale averaged $\Ld^\infty-\Ld^1$ estimate with a correction term involving the energy norm.
There is obviously a trade-off in $R$ in this estimate: the smaller~$R$, the stronger the left-hand side;
the larger~$R$, the smaller the correction term in the right-hand side. Note that this estimate~\eqref{eq:periodic-disp-0} still encodes dispersive effects in the regime of times~$t$ and averaging 
scales~$R$ for which the correction term can be absorbed into the left-hand side                                                                                                                                                                                                                                                                                                                             
(this is the way we shall proceed in the proof of Theorem~\ref{th:main}).
In view of the conservation of the $\Ld^2$ norm, this estimate therefore bears information if the energy norm is much smaller than the $\Ld^2$ norm, that is, at low Fourier modes.
This contrasts with the standard way how dispersive estimates such as~\eqref{eq:standard-disp-1} are exploited, that is, typically in form of Strichartz estimates to control high frequencies.

In order to take better advantage of~\eqref{eq:periodic-disp-0} at low frequencies, we note that it can be naturally post-processed by scaling.
More precisely, the above periodic argument shows that this estimate is unchanged upon rescaling the underlying periodic coefficient field.
{Applying~\eqref{eq:periodic-disp-0} with $u^\circ$, $t$, $R$, and the coefficient field~$a$ replaced by $u^\circ(\frac\cdot\kappa)$, $\kappa t$, $\kappa R$, and~$a(\frac\cdot\kappa)$, respectively, and changing variables, we obtain}
for all $t\ge0$, $\kappa,R>0$, and $u^\circ\in\Ld^1\cap H^1(\R^d)$,
\begin{multline}\label{eq:periodic-disp-0+}
\sup_{x\in\R^d}{R^{-\frac d2}}\|\!\cos(t\sqrt{H_a})u^\circ\|_{\Ld^2({B_{ R}}(x))}
\,\lesssim\,\kappa^{d}(1+\kappa t)^{-\frac{d-1}2}\|u^\circ\|_{\Ld^1(\R^d)}\\
+(1\wedge\lambda_0)^{-\frac12}\kappa^{-1}R^{-\frac d2}\|\nabla u^\circ\|_{\Ld^2(\R^d)}.
\end{multline}
The additional frequency parameter $\kappa$ in this estimate further increases the possible nontrivial trade-off between the terms.

We are now in position to state our main result on large-scale dispersive estimates. More precisely, we show that the estimate~\eqref{eq:periodic-disp-0+} indeed holds in the periodic setting (yet, in a slightly weaker version, with an additional error term $\kappa R^{-d/2}\|u^\circ\|_{\Ld^2(\R^d)}$ in the right-hand side), and we derive a corresponding improved estimate for localized data in the spirit of~\eqref{eq:standard-disp-2}. Although this could be proved by means of Floquet--Bloch theory as suggested above starting from~\eqref{eq:periodic-disp-Pc}, we provide a different  approach based on homogenization.
As it does not rely on Floquet--Bloch theory, our method can also be used to obtain corresponding results in the quasiperiodic and random settings.

\begin{theor}[Large-scale dispersive estimates]\label{th:disp}
Let $d\ge1$, consider a coefficient field {$A:\R^d\to\R^{d\times d}$} that is symmetric and uniformly elliptic in the sense of~\eqref{eq:ell-cond},
denote by~$U(t):=\cos(t\sqrt{H_a})$ the associated heterogeneous wave flow on $\Ld^2(\R^d)$, and 
consider an initial condition $u^\circ\in \Ld^1\cap H^1(\R^d)$.
Let also $\chi\in C^\infty(\R^d)$ be a smooth cut-off with $\chi(x)=1$ for $|x|<\frac12$ and $\chi(x)=0$ for $|x|>1$, and set $\chi_L:=\chi(\frac \cdot L)$ for $L\ge1$.
\begin{enumerate}[(i)]
\item \emph{Periodic setting:}\\ 
Assume that $A$ is periodic.
Then there exists a constant {$0<\kappa_*\le 1$} such that we have
for all $0<\kappa \le\kappa_*$ and $t,R\ge1$,
\begin{multline*}
\qquad\sup_{x\in\R^d} R^{-\frac d2}\|U(t)u^\circ\|_{\Ld^2({B_R}(x))}
\lesssim\,
 \kappa^d(1+\kappa t)^{-\frac{d-1}2}\|u^\circ\|_{\Ld^1(\R^d)}\\
+{R^{-\frac d2}}\Big(\kappa \|u^\circ\|_{\Ld^2(\R^d)}+\kappa^{-1}\|\nabla u^\circ\|_{\Ld^2(\R^d)}\Big).
\end{multline*}
In addition, given $\eta>0$, for all $0<\kappa \le\kappa_*$, $R,L\ge1$, and $t\gg R+L$,
\begin{multline*}
\qquad R^{-\frac d2}\|U(t)(\chi_Lu^\circ)\|_{\Ld^2({B_R})}\,
\lesssim_\eta \,
\kappa^d(1+\kappa t)^{\eta-d}\|u^\circ\|_{\Ld^1(B_L)}
\\+R^{-\frac d2}\Big(\big(\kappa+(\kappa L)^{-1}\big)\|u^\circ\|_{\Ld^2(B_L)}+\kappa^{-1}\|\nabla u^\circ\|_{\Ld^2(B_L)}\Big).
\end{multline*}

\item \emph{Quasiperiodic setting:}\\
Assume that $A$ is quasiperiodic with Diophantine and smoothness conditions as in the statement of Theorem~\ref{th:main}(ii).
Then there exist an exponent $\theta>0$ and a constant~$0<\kappa_* \le 1$ (only depending on the Diophantine and smoothness conditions on $A$), such that the following holds:
for all $0<\kappa \le\kappa_*$, $R\ge1$, and $ t \le \exp(\tfrac1C\kappa^{-\frac1{\theta}})$,
\begin{multline*}
\qquad \sup_{x\in\R^d}{R^{-\frac d2}} \|U(t)u^\circ\|_{\Ld^2({B_R}(x))}
\,\lesssim\,
\kappa ^d(1+\kappa t)^{-\frac{d-1}2}\|u^\circ\|_{\Ld^1(\R^d)}\\
+R^{-\frac d2}\Big(\kappa \|u^\circ\|_{\Ld^2(\R^d)}+\kappa^{-1}\|\nabla u^\circ\|_{\Ld^2(\R^d)}\Big).
\end{multline*}
In addition, given $\eta>0$, for all $0<\kappa \le\kappa_*$, $R,L\ge1$, and $R+L \ll t \le\exp(\tfrac1C\kappa^{-\frac1{\theta}})$,
\begin{multline*}
\qquad R^{-\frac d2} \|U(t)(\chi_Lu^\circ)\|_{\Ld^2(B_R)}
\,\lesssim_\eta\,
\kappa ^d(1+\kappa  t)^{\eta-d}\|u^\circ\|_{\Ld^1(B_L)}\\
+R^{-\frac d2}\Big(\big(\kappa+(\kappa L)^{-1}\big)\|u^\circ\|_{\Ld^2(B_L)}+\kappa^{-1}\|\nabla u^\circ\|_{\Ld^2(B_L)}\Big).
\end{multline*}

\item \emph{Random setting:}\\
Assume that $A=\tilde A(\cdot,\omega_\circ)$ is a realization of a random field $\tilde A$ as in the statement of Theorem~\ref{th:main}(iii),
and set 
\[{\alpha(d)\,:=\, \frac d2+1.}\]
For all $\e>0$, there exists a nonnegative, statistically spatially homogeneous random field $0<\kappa_{*,\e}\le 1$ with finite (inverse) stretched exponential moments,
\[\qquad\expecm{\exp(C^{-1}\kappa_{*,\e}^{-1/C})}\,<\,\infty,\qquad\text{for some $C<\infty$ (depending on $d,C_0,\e$)},\]
such that the following holds: for $\Pm$-almost all $\omega_\circ$, 
for all $x\in\R^d$, $0<\kappa \le  \kappa_{*,\e}(x)$, $0<r\le1$, $t\ge0$, and $R\ge1$ with $t,R \le\kappa^{r-\alpha(d)}$,
\begin{multline*}
\qquad R^{-\frac d2} \|U(t)u^\circ\|_{\Ld^2(B_R(x))}\,\lesssim\,
\kappa^d(1+\kappa t)^{-\frac{d-1}2}\|u^\circ\|_{\Ld^1(\R^d)}
\\
+R^{-\frac d2} \Big( \kappa_{*,\e}(x)^{-2}\kappa^{r-\e}\|u^\circ\|_{\Ld^2(\R^d)}+\kappa^{-1}\|\nabla u^\circ\|_{\Ld^2(\R^d)}\Big).
\end{multline*}
In addition, given $\eta>0$, for all $0<\kappa <\kappa_{*,\e}(0)$, $R,L\ge1$, and $R+L \ll t \le \kappa^{r-\alpha(d)}$,
\begin{align*}
&\quad R^{-\frac d2} \|U(t)(\chi_Lu^\circ)\|_{\Ld^2(B_R)}\,\lesssim_\eta\, 
\kappa^d(1+\kappa t)^{\eta-d} \|u^\circ\|_{\Ld^1(B_L)}\\
&\hspace{2.3cm}+R^{-\frac d2} \Big( \big(\kappa_{*,\e}(0)^{-2}\kappa^{r-\e}+(\kappa L)^{-1}\big)\|u^\circ\|_{\Ld^2(B_L)}+\kappa^{-1}\|\nabla u^\circ\|_{\Ld^2(B_L)}\Big).
\end{align*}
\end{enumerate}
\end{theor}

Note that quite some attention has been devoted in the literature to the validity of dispersive estimates for waves in heterogeneous media. We refer to~\cite{Schlag-07} and references therein for dispersive estimates for Schr\"odinger operators with a decaying potential field. Closer to the statistically spatially homogeneous setting considered in the present work, we also refer to~\cite{Bambusi-Zhao}, where dispersive estimates for 1D quasiperiodic Schr\"odinger operators are deduced from the very precise known description of the spectrum of such operators.
In these references, dispersive estimates are exclusively sought in the standard form~\eqref{eq:standard-disp-1}, possibly up to a projection on the continuous part of the spectrum as in~\eqref{eq:periodic-disp-Pc}. To our knowledge, the above result is the first investigation of weak dispersive estimates
with error terms involving the energy.

\subsection*{Structure of the article}
The rest of the article is organized as follows.
First, since it is not
primarily addressed
to experts in quantitative homogenization, we start in Section~\ref{sec:corr} by briefly recalling and motivating the needed input from the recent results on the long-time homogenization of the wave equation.
We then proceed in Section~\ref{sec:disp} to the proof of Theorem~\ref{th:disp} on large-scale dispersive estimates, and in Section~\ref{sec:th} to the proof of Theorem~\ref{th:main} on the spreading of potential eigenstates in the lower spectrum.

\section{A primer on long-time homogenization for waves}\label{sec:corr}
\begingroup\allowdisplaybreaks
Homogenization theory aims at describing the behavior of solutions of PDEs with heterogeneous coefficients in the regime when the typical spatial scale of the coefficients is much smaller than the spatial scale of external forces and initial data: in such a regime, one expects the coefficients to average out in some sense, so that solutions of the PDE would be described to leading order by solutions of a corresponding `homogenized' constant-coefficient PDE (possibly of a different type). We refer
to~\cite{BLP-78} for a systematic treatment of homogenization of various PDEs with periodic coefficients.

For the purpose of the present work, we shall consider the wave equation associated with the heterogeneous acoustic operator $H_a=-\nabla\cdot A\nabla$,
\begin{equation}\label{eq:0-wave}
\left\{\begin{array}{l}
\partial_{t}^2u_\kappa=\nabla\cdot A\nabla u_\kappa,\\
u_\kappa|_{t=0}=u^\circ_\kappa,\\
\partial_tu_\kappa|_{t=0}=0,
\end{array}\right.
\end{equation}
with some (smooth) initial condition~$u^\circ_\kappa \in \Ld^2(\R^d)$.
If the initial condition is localized on low frequencies $0<\kappa\ll1$ (say for instance $u_\kappa^\circ=\kappa^{d/2} u^\circ(\kappa\cdot)$), then there is a scale separation with respect to variations of~$A$ on the unit scale,
so that homogenization is expected to hold --- under suitable statistical spatial homogeneity assumptions on the coefficient field~$A$ (for instance, if $A$ is periodic, quasiperiodic, or if~$A$ is a typical realization of a random field with statistical spatial homogeneity as in Theorem~\ref{th:main}(iii)).
More precisely, as shown e.g.~in~\cite{Brahim-Otsmane-92} in the periodic setting, it is a standard result of homogenization theory that~$u_\kappa$ is then close to the solution~$\bar u_\kappa$ of the homogenized wave equation
\begin{equation}\label{eq:0-wave-homg}
\left\{\begin{array}{l}
\partial_t^2\bar u_\kappa=\nabla\cdot \bar A\nabla \bar u_\kappa,\\
\bar u_\kappa |_{t=0}=u^\circ_\kappa,\\
\partial_t\bar u_\kappa|_{t=0}=0,
\end{array}\right.
\end{equation}
{on any compact time interval $[0,\kappa^{-1}T]$,} where the homogenized coefficient~$\bar A$ is some constant matrix.
Recall that $\bar A$ only depends on $A$ and is in particular independent of initial data. What is less standard and was only established recently, is the quantification of the maximal timescale up to which $u_\kappa$ indeed remains close to $\bar u_\kappa$. In the periodic setting, as first understood in~\cite{SS-91}, this holds only {up to times~$t\ll \kappa^{-3}$,} while on longer timescales the homogenized equation~\eqref{eq:0-wave-homg} is no longer accurate and dispersive corrections need to be added to the homogenized operator $-\nabla\cdot \bar A\nabla$.
More precisely, as shown in~\cite{BG,DGR}, there exist higher-order constant tensors $\{\bar A^{n}\}_{n\ge1}$ with $\bar A^1=\bar A$ such that for all \mbox{$N\in\N$} the heterogeneous solution $u_\kappa$ is well-approximated up to times $t\ll\kappa^{-N-1}$ by the solution~$\bar u_{\kappa,N}$ of the corrected homogenized wave equation
\begin{equation}\label{eq:0homog-baru}
\left\{\begin{array}{l}
\partial_t^2\bar u_{\kappa,N}=\nabla\cdot\big(\sum_{n=1}^N\bar A^n_{j_1\ldots j_{n-1}}\nabla^{n-1}_{j_1\ldots j_{n-1}}\big)\nabla\bar u_{\kappa,N},\\
\bar u_{\kappa,N}|_{t=0}=u^\circ_{\kappa},\\
\partial_t \bar u_{\kappa,N}|_{t=0}=0.
\end{array}\right.
\end{equation}
Of course, the maximal timescale up to which the result holds strongly depends on assumptions on $A$: whereas in the periodic setting one can take any~$N\in\N$, the result can only hold up to some maximal exponent $N_0$ in the random setting (which in 1D can be related to Anderson localization~\cite{Sims-Stolz}).

In order to motivate the upcoming results, let us quickly describe informally why~$u_\kappa$ might indeed be described by the solution~$\bar u_{\kappa,N}$ of a homogenized equation. To this aim, we follow the approach introduced in~\cite{BG} and further developed in~\cite{DGR}, based on so-called ``approximate spectral theory'' {(that builds upon the first results~\cite{MR2843023,MR3191584} in the field)}.
Starting point is to try diagonalizing the acoustic operator~$H_a=-\nabla \cdot A \nabla$ at the bottom of its spectrum by means of extended states in form of Bloch waves $x \mapsto e^{ix\cdot\xi} \chi_\xi(x)$, where~$\chi_\xi$ would have the same type of statistical spatial homogeneity as $A$.
In the periodic case, this is precisely realized by Floquet theory, see e.g.~\cite{Kuchment-16}, whereas in the random case statistically homogeneous Bloch waves are not expected to exist~\cite{DS21}. With this failure of Floquet theory, we are led to rather look for an approximate diagonalization of $H_a$.
More precisely, we start with the following observation: in the periodic setting, as shown e.g.\@ in~\cite{ABV-16}, the gradient of the Bloch wave~$\chi_\xi$ with respect to the wavenumber $\xi$ at $0$ coincides with the so-called corrector $\phi$ from elliptic homogenization theory, that is, the periodic solution (with vanishing average) of 
\[-\nabla \cdot A(\nabla \phi+\Id)=0.\]
 In other words, one has a Taylor expansion $\chi_\xi=1+i\xi_j\phi_j+O(|\xi|^2)$ for the Bloch wave. 
As shown in~\cite{BG}, this expansion can be pursued to higher orders and we find
\[\chi_\xi\,=\,\sum_{n=0}^N (i\xi)^{\otimes n}_{j_1\ldots j_n}\phi^n_{j_1\ldots j_n}+O(|\xi|^{N+1})\]
for all $N\in\N$ in the periodic setting, where the coefficients $\{\phi^n\}_n$ are so-called higher-order (spectral) correctors and satisfy some hierarchy of PDEs.
The main idea in~\cite{BG} was then to replace the use of the Bloch wave $\chi_\xi$ by its higher-order expansions for all practical purposes: the advantage is that
 correctors $\{\phi^n\}_n$ can be defined in settings when Bloch waves do not exist (or are not known to exist) --- for instance, all of them can easily be defined in the quasiperiodic setting (under some Diophantine condition), and at least a finite number of them in the random setting (under suitable mixing assumptions).
The drawback is of course that these expansions do not exactly diagonalize the acoustic operator $H_a=-\nabla \cdot A \nabla$: there is a remainder in the generalized eigenvalue equation (called eigendefect in~\cite{BG}). Still, this yields an approximate diagonalization of~$H_a$ and we are led to an approximation of~$u_\kappa$ by means of a so-called two-scale expansion of the form
\begin{equation}\label{eq:expand-2scale}
\sum_{n=0}^N \phi^n_{j_1\ldots j_n}\nabla^n_{j_1\ldots j_n}\bar u_{\kappa,N},
\end{equation}
for some (smooth) effective profile $\bar u_{\kappa,N}$. Inserting this ansatz into~\eqref{eq:0-wave} easily shows that this profile~$\bar u_{\kappa,N}$ must solve an homogenized wave equation of the form~\eqref{eq:0homog-baru}. For initial data~$u_\kappa^\circ$ localized on low frequencies $0<\kappa\ll1$,
we find $\nabla^n\bar u_{\kappa,N}=O(\kappa^n)$ in $\Ld^2(\R^d)$, so the above expansion~\eqref{eq:expand-2scale} coincides with~$\bar u_{\kappa,N}$ to leading order. This formally justifies why~$\bar u_{\kappa,N}$ indeed provides a good description of~$u_\kappa$ in $\Ld^2(\R^d)$. Yet, the whole approach only works provided that we can control the error made by using an approximate diagonalization of $H_a$, and this depends on the size of the remainder in the eigenvalue equation: the smaller this error, the longer the time the approximate solution~$\bar u_{\kappa,N}$ remains close to~$u_\kappa$. 
If the corrector $\phi^{N}$ is a bounded function, the remainder in the eigenvalue equation is essentially $O(\kappa^{N+1})$, and then using energy estimates on the wave equation leads us to an approximation error $u_\kappa-\bar u_{\kappa,N}=O(\kappa+t \kappa^{N+1})$ in $\Ld^2(\R^d)$.
This is the origin for the timescale condition~$t\ll\kappa^{-N-1}$ mentioned above. 

The quantitative analysis of corrector equations has been central in modern homogenization, and we have by now very accurate bounds on correctors $\{\phi^n\}_n$ in various settings, building up on~\cite{MR512007} in the quasiperiodic setting and on~\cite{GO1,Gloria-Otto-10b,AKM2} in the random setting.
Before we state these bounds, let us give the precise form of the hierarchy of corrector equations.
To cover all settings at once, we formulate these equations in the whole space~$\R^d$, while implicitly understanding that we look for periodic solutions in the periodic setting, for quasiperiodic solutions in the quasiperiodic setting, and for statistically spatially homogeneous solutions (when possible) in the random setting. We also denote by $\expec{\cdot}$ the `mean' or infinite-volume average, which is understood as average over the periodicity cell in the periodic setting, average over the underlying high-dimensional periodicity cell in the quasiperiodic setting, and expectation in the random setting.
With this notation, correctors $\phi^n=(\phi^n_{j_1\ldots j_n})_{1\le j_1,\ldots,j_n\le d}\in H^1_\loc(\R^d;(\R^{d})^n)$ are defined inductively by~$\phi^0=1$ and by the following equations, for all $n\ge1$,
\begin{multline}\label{eq:correctors}
\qquad-\nabla\cdot A\nabla\phi^n_{j_1\ldots j_n}\,=\,\nabla\cdot \big(A\phi^{n-1}_{j_1\ldots j_{n-1}}e_{j_{n}}\big)+e_{j_{n}}\cdot A\big(\nabla\phi^{n-1}_{j_1\ldots j_{n-1}}+\phi^{n-2}_{j_1\ldots j_{n-2}}e_{j_{n-1}}\big)\\
-\sum_{m=1}^{n-1}e_{j_n}\cdot\bar A^m_{j_1\ldots j_{m-1}}\phi^{n-m-1}_{j_m\ldots j_{n-2}}e_{j_{n-1}},
\end{multline}
where we choose $\phi^n_{j_1\ldots j_n}$ with zero mean $\E[\phi^n_{j_1\ldots j_n}]=0$ (when possible),
and where the constant tensors $\{\bar A^m\}_m$ are iteratively chosen to make sure that the right-hand side of~\eqref{eq:correctors} also has zero mean. More precisely, this amounts to defining for all $m\ge1$,
\[\bar A^m_{j_1\ldots j_{m-1}}e_{j_m}\,:=\,\expecm{A\big(\nabla \phi^{m}_{j_1\ldots j_m}+\phi_{j_1\ldots j_{m-1}}^{m-1}e_{j_m}\big)}.\]
We implicitly set $\phi^{-1}=0$ for notational convenience.
As usual in homogenization theory, it is convenient to further introduce flux correctors $\sigma^n:=(\sigma^n_{j_1\ldots j_n})_{1\le j_1,\ldots,j_n\le d}\in H^1_\loc(\R^d;\R^{d\times d}\times(\R^d)^n)$:
for all~$n\ge1$, we define
\[\sigma^n_{j_1\ldots j_n}=\nabla \Phi_{j_1\ldots j_n}^n,\]
where $\Phi_{j_1\ldots j_n}^n\in H^2_\loc(\R^d;\R^d)$ satisfies
\[\triangle\Phi_{j_1\ldots j_n}^n=A\big(\nabla\phi^{n}_{j_1\ldots j_{n}}+\phi^{n-1}_{j_1\ldots j_{n-1}}e_{j_{n}}\big)-\sum_{m=1}^n\bar A^m_{j_1\ldots j_{m-1}}\phi^{n-m}_{j_m\ldots j_{n-1}}e_{j_{n}}.\]
In particular, if well-defined, we note that
this construction ensures for all $n\ge1$,
\begin{equation}\label{eq:def-sigma}
\nabla\cdot\sigma_{j_1\ldots j_n}^n=A\big(\nabla\phi^{n}_{j_1\ldots j_{n}}+\phi^{n-1}_{j_1\ldots j_{n-1}}e_{j_{n}}\big)-\sum_{m=1}^n\bar A^m_{j_1\ldots j_{m-1}}\phi^{n-m}_{j_m\ldots j_{n-1}}e_{j_{n}}.
\end{equation}
We also recall that $\bar A^n=0$ for $n$ even, if well-defined, cf.~\cite{BG,DGR}, which ensures the symmetry of the differential operator in the homogenized wave equation~\eqref{eq:0homog-baru}.
Quantitative homogenization theory provides the following optimal corrector estimates; note how the result strongly varies between the periodic, quasiperiodic, and random settings.
\begin{lem}[Correctors in periodic setting; e.g.~\cite{BLP-78,DGL-19}]\label{lem:cor-per}
Assume that the coefficient field $A:\R^d\to\R^{d\times d}$ is symmetric and uniformly elliptic in the sense of~\eqref{eq:ell-cond}, and is periodic on~$Q=[0,1)^d$. Then, for all $n\in\N$, the correctors $(\phi^n,\sigma^n)$ are uniquely defined as mean-zero $Q$-periodic functions and satisfy
\[ \|(\phi^n,\sigma^n)\|_{H^1(Q)}\,\le\, C^n, \qquad |\bar A^n|\,\le\, C^n.\]
\end{lem}

\begin{lem}[Correctors in quasiperiodic setting; e.g.~\cite{MR512007,MR3451428,DGS-21}]\label{lem:cor-qper}
Assume that the coefficient field $A:\R^d\to\R^{d\times d}$ is symmetric and uniformly elliptic in the sense of~\eqref{eq:ell-cond}, and takes the form
$A(x)=A_0(Fx)$ for some frequency matrix $F\in\R^{M\times d}$ with $M>d$, and for some lifted map $A_0:\R^M\to\R^{d\times d}$ that is periodic on $\R^M$ and Gevrey-regular. Further assume that the frequency matrix $F$ satisfies a Diophantine condition, that is, for some~$C,\gamma>0$,
\[|Fz|\,\ge\,\tfrac1C|z|^{-\gamma},\qquad\text{for all $z\in\Z^M\setminus\{0\}$}.\]
Then, for all $n\in\N$, the correctors $(\phi^n,\sigma^n)$ are uniquely defined as mean-zero quasiperiodic functions and satisfy
\[ \|(\phi^n,\psi^n)\|_{W^{1,\infty}(\R^d)}\,\le\,(C n^\theta)^n, \qquad |\bar A^n|\,\le\, (C n^\theta)^n,\]
for some exponent $\theta>0$ only depending on $\gamma$ and on the Gevrey regularity of $A_0$.
\end{lem}

\begin{lem}[Correctors in random setting; e.g.~\cite{
GNO-quant,AKM-book,MR4147699,
DGR}]\label{lem:cor-rand}
Assume that the coefficient field is given by a symmetric, uniformly elliptic, statistically spatially homogeneous random field $A:\R^d\times\Omega\to\R^{d\times d}$, constructed on some probability space $(\Omega,\Pm)$, and further assume that it satisfies strong enough mixing assumptions, as for instance having finite range of dependence, in the sense of Theorem~\ref{th:main}(iii).
Then, the corrector gradients $(\nabla \phi^n,\nabla \sigma^n)$ are uniquely defined as mean-zero statistically spatially homogeneous random fields for all $1\le n\le \lceil \frac d2 \rceil$, such that the correctors $(\phi^n,\sigma^n)$ are themselves mean-zero statistically homogeneous fields for $1\le n< \lceil \frac d2 \rceil$, and such that corrector equations with coefficients $A(\cdot,\omega)$ are satisfied by $\{\phi^n(\cdot,\omega),\sigma^n(\cdot,\omega)\}_n$ for $\Pm$-almost all $\omega$.
The homogenized tensors $\bar A^m$ are then well-defined for all $1\le m\le\lceil \frac d2 \rceil$.
Moreover, for all $\e>0$ and $\delta>\frac12$, there exists a random variable $\Kc_{\e,\delta}$ with finite stretched-exponential moments,
\[\expecm{\exp(C^{-1}\Kc_{\e,\delta}^{1/C})}<\infty,\qquad\text{for some $C<\infty$ (depending on $d,C_0,\e,\delta$)},\]
such that, $\Pm$-almost surely,
\begin{enumerate}[---]
\item for all $n\le\lfloor \frac d2\rfloor$,
\[\quad|(\phi^n,\sigma^n)(x)|+|(\nabla \phi^n,\nabla \sigma^n)(x)|\,\le\, \Kc_{\e,\delta}\langle x\rangle^\e;\]
\item for $d$ odd and $n= \lceil \frac d2\rceil$,
\begin{equation}\label{e.wo-massive-app}
\quad|(\phi^n,\sigma^n)(x)|\,\le\, \Kc_{\e,\delta}\langle x\rangle^\delta\qquad
\text{and}\qquad|(\nabla \phi^n,\nabla \sigma^n)(x)|\,\le\, \Kc_{\e,\delta}\langle x\rangle^\e.
\end{equation}
\end{enumerate}
\end{lem}
\medskip

{For the proof of long-time homogenization, and in particular in our application to large-scale dispersive estimates in the random setting, the algebraic growth of the corrector of order $n= \lceil \frac d2\rceil$ (when~$d$ is odd) is problematic and only leads to pessimistic error bounds. Instead, it is more advantageous to introduce an artificial cut-off at low frequencies (e.g.\@ in the spirit of the approximation of homogenized coefficients in~\cite{MR3451428}).
Such a cut-off was not used in previous works on long-time homogenization~\cite{BG,DGR} and was inspired to us by a recent related work on the long-wave approximation of random FPUT lattices~\cite{McGinnis-Wright-23}; a detailed account of this improvement is postponed to a separate note, which we briefly summarize here.
Next to the above-defined corrector $\phi^n$, for $n\le\lceil \frac d2\rceil$, in the random setting of Lemma~\ref{lem:cor-rand}, we define its massive version $\phi^n_\kappa$ as the unique statistically spatially homogeneous solution of 
\begin{multline}\label{eq:correctors-mass}
\kappa^2\phi^n_{\kappa,j_1\ldots j_n} -\nabla\cdot A\nabla\phi^n_{\kappa,j_1\ldots j_n}\,=\,\nabla\cdot \big(A\phi^{n-1}_{j_1\ldots j_{n-1}}e_{j_{n}}\big)+e_{j_{n}}\cdot A\big(\nabla\phi^{n-1}_{j_1\ldots j_{n-1}}+\phi^{n-2}_{j_1\ldots j_{n-2}}e_{j_{n-1}}\big)\\
-\sum_{m=1}^{n-1}e_{j_n}\cdot\bar A^m_{j_1\ldots j_{m-1}}\phi^{n-m-1}_{j_m\ldots j_{n-2}}e_{j_{n-1}},
\end{multline}
which corresponds to the corrector equation~\eqref{eq:correctors} with a massive term $\kappa^2$.
Note that the well-posedness of this equation is ensured by the massive term.
This gives rise to the massive approximation $\bar A^n_\kappa$ of $\bar A^n$ that we define by
\[
\bar A^n_{\kappa,j_1\ldots j_{n-1}}e_{j_n}\,:=\,\expecm{A\big(\nabla \phi^{n}_{\kappa,j_1\ldots j_n}+\phi_{j_1\ldots j_{n-1}}^{n-1}e_{j_n}\big)}.
\]
Because of the massive term in \eqref{eq:correctors-mass}, the flux $A\big(\nabla \phi^{n}_{\kappa,j_1\ldots j_n}+\phi_{j_1\ldots j_{n-1}}^{n-1}e_{j_n}\big)$ is not divergence-free, and one cannot put 
\[A\big(\nabla \phi^{n}_{\kappa,j_1\ldots j_n}+\phi_{j_1\ldots j_{n-1}}^{n-1}e_{j_n}\big)-\bar A^n_{\kappa,j_1\ldots j_{n-1}}e_{j_{n}}-\sum_{m=1}^{n-1}\bar A^m_{j_1\ldots j_{m-1}}\phi^{n-m}_{j_m\ldots j_{n-1}}e_{j_{n}}\]
exactly in divergence form as we do in~\eqref{eq:def-sigma}. Similarly as in~\cite{MR4103433}, there is a defect and we are led to the decomposition 
\begin{multline}\label{e.new-flux}
A\big(\nabla \phi^{n}_{\kappa,j_1\ldots j_n}+\phi_{j_1\ldots j_{n-1}}^{n-1}e_{j_n}\big)-\bar A^n_{\kappa,j_1\ldots j_{n-1}}e_{j_{n}}-\sum_{m=1}^{n-1}\bar A^m_{j_1\ldots j_{m-1}}\phi^{n-m}_{j_m\ldots j_{n-1}}e_{j_{n}}
\\
=\nabla \cdot \sigma^n_{\kappa,j_1\ldots j_n} + \kappa g^n_{\kappa,j_1\ldots j_n} ,
\end{multline}
where we define
\[\sigma^n_{\kappa,j_1\ldots j_n}:=\nabla \Phi^n_{\kappa,j_1\ldots j_n}\qquad\text{and}\qquad g_{\kappa,j_1\ldots j_n}^n=-\kappa\Phi^n_{\kappa,j_1\ldots j_n},\]
in terms of the unique statistically spatially homogeneous solution $\Phi_{\kappa,j_1\ldots j_n}$ of
\begin{multline*}
\kappa^2 \Phi_{\kappa,j_1\ldots j_n}^n-\triangle\Phi_{\kappa,j_1\ldots j_n}^n
\\
=\, -A\big(\nabla \phi^{n}_{\kappa,j_1\ldots j_n}+\phi_{j_1\ldots j_{n-1}}^{n-1}e_{j_n}\big)+\bar A^n_{\kappa,j_1\ldots j_{n-1}}e_{j_{n}}+\sum_{m=1}^{n-1}\bar A^m_{j_1\ldots j_{m-1}}\phi^{n-m}_{j_m\ldots j_{n-1}}e_{j_{n}}.
\end{multline*}
In the setting of Lemma~\ref{lem:cor-rand}, for $n=\lceil\frac d2\rceil$, we have, instead of~\eqref{e.wo-massive-app},
\begin{eqnarray}
|(\phi^n_\kappa,\sigma^n_\kappa, g^n_\kappa)(x)|&\le& \Kc_{\e,\delta}\langle x\rangle^\e \kappa^{-\delta},\label{e.massive-app}\\
|(\nabla \phi^n_\kappa,\nabla \sigma^n_\kappa)(x)|&\le&\Kc_{\e,\delta}\langle x\rangle^\e,\nonumber\\
 |\bar A^n_\kappa-\bar A^n|&\lesssim&\kappa,\nonumber
\end{eqnarray}
where now $\delta=0$ if $d$ is even and $\delta=\frac12$ if $d$ is odd.
In the periodic and quasiperiodic settings, the massive correctors remain uniformly bounded.}

\section{Proof of the large-scale dispersive estimates}\label{sec:disp}
This section is devoted to the proof of Theorem~\ref{th:disp}.
Given an initial condition $u^\circ\in \Ld^1\cap H^1(\R^d)$ with normalization $\|u^\circ\|_{\Ld^2(\R^d)}=1$, consider the solution $u^t=U(t)u^\circ$ of the heterogeneous wave equation
\begin{equation}\label{eq:wave-u0}
\left\{\begin{array}{l}
\partial_t^2u=\nabla\cdot A\nabla u,\\
u|_{t=0}=u^\circ,\\
\partial_tu|_{t=0}=0.
\end{array}\right.
\end{equation}
We split the proof into six main steps.

\medskip
\step1 Regularization of initial data.\\
We wish to appeal to homogenization theory for~\eqref{eq:wave-u0} while only assuming that the initial data~$u^\circ$ has low energy $\|\nabla u^\circ\|_{\Ld^2(\R^d)}\ll1$. Long-time homogenization results~\cite{BG,ALR,DGR} however require a finer control on the low-frequency nature of initial data, typically in form of improved smallness for higher derivatives $\nabla^k u^\circ$, $k\ge1$.
To deal with this issue, we start by defining a suitable regularization of initial data $u^\circ$:
given a nonnegative cut-off function $\chi\in C^\infty_c(B)$ with $\chi=1$ on $\frac12B$, and given a nonnegative convolution kernel $\rho\in\Sc(\R^d)$ with 
{$\int_{\R^d}\rho=1$}, whose Fourier transform $\hat\rho$ is supported in $B$, 
we consider the rescalings
\[\chi_L\,:=\,\chi(\tfrac\cdot L),\qquad
\rho_\kappa\,:=\,\kappa^{d}\rho(\kappa\cdot),\]
where
$L\ge1$ and
{$0<\kappa \le 1$} will be properly chosen later on in the proof. We then define the regularized initial condition
\begin{equation}\label{eq:def-psilambda}
u_{\kappa,L}^\circ\,:=\,\rho_\kappa\ast (\chi_Lu^\circ),
\end{equation}
and we consider the corresponding wave equation,
\begin{equation}\label{eq:wave-u0-reg}
\left\{\begin{array}{l}
\partial_t^2u_{\kappa,L}=\nabla\cdot A\nabla u_{\kappa,L},\\
u_{\kappa,L}|_{t=0}=u_{\kappa,L}^\circ,\\
\partial_tu_{\kappa,L}|_{t=0}=0.
\end{array}\right.
\end{equation}
Comparing with~\eqref{eq:wave-u0} and appealing to the finite speed of propagation $c \simeq 1$, {which implies that at time $t$ we have $(U(t) u^\circ) |_{B_R} = (U(t) (u^\circ \mathds1_{B_{R+ct}}) |_{B_R}$,} we find
for all~$R,t\ge0$ and $L \ge 4(R+c t)$,
\begin{eqnarray*}
\|u^t-u_{\kappa,L}^t\|_{\Ld^2(B_R)}
&\le&\|u^\circ-u_{\kappa,L}^\circ\|_{\Ld^2(B_{R+ct})}\\
&\le&\|\rho_\kappa*(u^\circ(1-\chi_L))\|_{\Ld^2(B_{L/4})}
+\|u^\circ-\rho_\kappa\ast u^\circ\|_{\Ld^2(\R^d)}.
\end{eqnarray*}
{Using that $1-\chi_L$ vanishes on $B_{L/2}$, using Jensen's inequality with $\int_{\R^d}\rho_\kappa=1$, and using that $\rho$ has superalgebraic decay at infinity, we find for the first right-hand side term, for any $n\ge1$,
\begin{eqnarray*}
\|\rho_\kappa*(u^\circ(1-\chi_L))\|_{\Ld^2(B_{L/4})}^2&=&\int_{B_{L/4}} \Big(\int_{\R^d} \rho_\kappa(x-y)u^\circ(y)(1-\chi_L(y))\,dy\Big)^2 dx\\
&\le& \int_{B_{L/4}} \Big(\int_{\R^d \setminus B_{L/2}}  \rho_\kappa(x-y)  |u^\circ(y)|\, dy\Big)^2 dx
\\
&\le& \Big(\int_{\R^d} |u^\circ|^2\Big) \int_{\R^d \setminus B_{L/4}}\rho_\kappa(y)\,dy\\
&\lesssim_n& (\kappa L)^{-2n} \|u^\circ\|_{\Ld^2(\R^d)}^2.
\end{eqnarray*}
Further noting that
\[\|u^\circ-\rho_\kappa\ast u^\circ\|_{\Ld^2(\R^d)}\,\le\,\Big(\int_{\R^d}|y|^2\rho_\kappa(y)\,dy\Big)^\frac12\|\nabla u^\circ\|_{\Ld^2(\R^d)}\,\lesssim\,\kappa^{-1}\|\nabla u^\circ\|_{\Ld^2(\R^d)},\]
the above becomes, for all $R,t\ge0$, $L \ge 4(R+c t)$, and $n\ge1$,}
\begin{equation}\label{eq:err-init-reg}
\|u^t-u_{\kappa,L}^t\|_{\Ld^2(B_R)}
\,\le\,C_n(\kappa L)^{-n}+C\kappa^{-1}\|\nabla u^\circ\|_{\Ld^2(\R^d)}.
\end{equation}
If the initial condition $u^\circ$ is supported in $B_{L/2}$,
we do not need the finite speed of propagation and we simply find for all $t,R\ge 0$,
\begin{equation}
\|u^t-u_{\kappa,L}^t\|_{\Ld^2(B_R)}
\,\le\,C \kappa^{-1}\|\nabla u^\circ\|_{\Ld^2(\R^d)}.\label{eq:err-init-reg+}
\end{equation}

\medskip
\step2 Higher-order correctors and homogenized equations.\\
We assume that we can iteratively construct $N\ge1$ controlled correctors {$\{\phi^n,\sigma^n\}_{1\le n< N}$ 
and $\{\phi^N_\kappa,\sigma^N_\kappa,g^N_\kappa\}$,} and homogenized tensors {$\{\bar A^n\}_{1\le n<N}$ and $\bar A^N_\kappa$} in the sense of Section~\ref{sec:corr} with the following pointwise bounds,
\begin{gather}
|\bar A^n|\le K_n \quad\text{for all~\mbox{$1\le n<N$}}, \quad {|\bar A^N_\kappa|\le K_N,} \label{eq:estim-corr}\\
|(\phi^n,\sigma^n)(x)|+|(\nabla\phi^{n},\nabla\sigma^{n})(x)|\le K_n\langle x\rangle^\e,\quad\text{for all~\mbox{$1\le n<N$}},\nonumber\\
{|(\phi^N_\kappa,\sigma^N_\kappa, g_\kappa^N)(x)|\le K_N\langle x\rangle^\e \kappa^{-\delta},\quad |(\nabla\phi^{N}_\kappa,\nabla\sigma^{N}_\kappa)(x)|\le K_n\langle x\rangle^\e} \nonumber
\end{gather}
for some constants~$K_0,\ldots,K_N\ge1$ and some exponents~\mbox{$0\le \e\le\delta<1$}.
In the sequel of the proof, we argue by assuming that those properties are satisfied in this general form, before particularizing them to the specific case of periodic, quasiperiodic, and random coefficient fields, as described in Section~\ref{sec:corr}. In a nutshell:
\begin{enumerate}[---]
\item in the periodic case, the above bounds hold with $N\uparrow\infty$, $\e=\delta=0$, and~$K_n=C^n$;
\smallskip\item in the quasiperiodic case, under a Diophantine condition and Gevrey regularity, the above bounds hold with $N\uparrow\infty$, $\e=\delta=0$, and $K_n=(Cn^\theta)^n$ for some $\theta>0$;
\smallskip\item in the random case, under suitable mixing assumptions, the above bounds hold with only $N=\lceil\frac d2\rceil$, with any $\e>0$, {and with $\delta=0$ if the spatial dimension $d$ is even and~$\delta=\frac12$ if $d$ is odd.}
\end{enumerate}
{For notational convenience, when there is no confusion, we drop the subscript $\kappa$ in correctors and homogenized tensors of order $N$.}

\medskip\noindent
In terms of the above homogenized tensors, we consider the associated (constant-coefficient) homogenized wave equation
\begin{equation}\label{eq:homog-baru}
\left\{\begin{array}{l}
\partial_t^2\bar u_{\kappa,L,N}=\nabla\cdot\big(\sum_{n=1}^N\bar A^n_{j_1\ldots j_{n-1}}\nabla^{n-1}_{j_1\ldots j_{n-1}}\big)\nabla\bar u_{\kappa,L,N},\\
\bar u_{\kappa,L,N}|_{t=0}=u_{\kappa,L}^\circ,\\
\partial_t \bar u_{\kappa,L,N}|_{t=0}=0.
\end{array}\right.
\end{equation}
This equation is however {ill-posed for general initial conditions} due to the lack of definiteness of the symbol
\begin{equation}\label{eq:def-muN}
\mu_N(i\xi):=\xi\cdot\Big(\sum_{n=1}^N(i\xi)^{\otimes(n-1)}_{j_1\ldots j_{n-1}}\bar A_{j_1\ldots j_{n-1}}^n\Big)\xi.
\end{equation}
To circumvent this issue, first recall that the uniform ellipticity condition~\eqref{eq:ell-cond} entails that $\bar A=\bar A^1$ also satisfies
\begin{equation}\label{eq:ell-barA}
\tfrac1{C_0}|e|^2\,\le\,e\cdot\bar Ae\,\le\,C_0|e|^2,\qquad\text{for all $e\in\R^d$},
\end{equation}
see e.g.~\cite{BLP-78}.
Therefore, using the above bounds~\eqref{eq:estim-corr},
we find for the symbol~$\mu_N$, for all~$\xi\in\kappa B$,
\begin{equation}\label{eq:lower-symb}
\mu_N(i\xi)\,\left\{\begin{array}{l}
\ge\,\xi\cdot\bar A\xi\big(1-C_0\sum_{n=2}^NK_n\kappa^{n-1}\big)\,\ge\,\tfrac12\xi\cdot\bar A\xi\,\ge\,\tfrac1{2C_0}|\xi|^2,\\
\vspace{-0.3cm}\\
\le\,\xi\cdot\bar A\xi\big(1+C_0\sum_{n=2}^NK_n\kappa^{n-1}\big)\,\le\,\tfrac32\xi\cdot\bar A\xi\,\le\,\tfrac{3C_0}{2}|\xi|^2,
\end{array}\right.
\end{equation}
provided that $\kappa$ is small enough in the sense of
\begin{equation}\label{eq:cond-small-lambda}
\sum_{n=2}^NK_n\kappa^{n-1}\,\le\,\tfrac1{2C_0}.
\end{equation}
(The symbol $\mu_N$ is obviously real-valued as $\bar A^n=0$ for $n$ even, cf.~\cite{BG,DGR}.)
Noting that by definition the initial Fourier transform~$\hat u_{\kappa,L}^\circ$ is supported in $\kappa B$, just as~$\hat \rho_\kappa$,  and using that by~\eqref{eq:lower-symb} the symbol is definite and bounded in $\kappa B$, equation~\eqref{eq:homog-baru} admits a unique solution given by the Fourier formula
\begin{equation}\label{eq:baru-lam}
\bar u_{\kappa,L,N}^t(x)\,:=\,\int_{\R^d}e^{ix\cdot\xi}\cos\big(t\sqrt{\mu_{N}(i\xi)}\big)\,\hat u^\circ_{\kappa,L}(\xi)\,d^*\xi.
\end{equation}

\medskip
\step3 Two-scale analysis for $u_{\kappa,L}$.\\
In terms of the above correctors, in the spirit of~\eqref{eq:expand-2scale}, we claim that the two-scale expansion
\[\sum_{n=0}^N\phi^n_{j_1\ldots j_n}\nabla^n_{j_1\ldots j_n}\bar u_{\kappa,L,N}\]
provides a good approximation for $u_{\kappa,L}$.
To prove that, we proceed by examining the wave equation satisfied by the error. 
Applying the wave operator, we obtain
\begin{multline*}
\big(\partial_t^2-\nabla\cdot A\nabla\big)\Big(u_{\kappa,L}-\sum_{n=0}^N\phi^n_{j_1\ldots j_n}\nabla^n_{j_1\ldots j_n}\bar u_{\kappa,L,N}\Big)
\,=\,
-\sum_{n=0}^N\phi^n_{j_1\ldots j_n}\nabla^n_{j_1\ldots j_n}\partial_t^2\bar u_{\kappa,L,N}\\
+\sum_{n=0}^N\nabla\cdot (A\nabla\phi^n_{j_1\ldots j_n})\nabla^n_{j_1\ldots j_n}\bar u_{\kappa,L,N}
+\sum_{n=0}^Ne_{j_{n+1}}\cdot A\nabla\phi^n_{j_1\ldots j_n}\nabla^{n+1}_{j_1\ldots j_{n+1}}\bar u_{\kappa,L,N}\\
+ \sum_{n=0}^N\nabla\cdot (A\phi^n_{j_1\ldots j_n}e_{j_{n+1}})\nabla^{n+1}_{j_1\ldots j_{n+1}}\bar u_{\kappa,L,N}
+ \sum_{n=0}^Ne_{j_{n+2}}\cdot A\phi^n_{j_1\ldots j_n}e_{j_{n+1}}\nabla^{n+2}_{j_1\ldots j_{n+2}}\bar u_{\kappa,L,N}.
\end{multline*}
Inserting the homogenized wave equation~\eqref{eq:homog-baru} for $\bar u_{\kappa,L,N}$ in the first right-hand side term for~$n<N$, relabeling the different sums, and recognizing the corrector equations~\eqref{eq:correctors} {and~\eqref{eq:correctors-mass}, this becomes
\begin{multline*}
\big(\partial_t^2-\nabla\cdot A\nabla\big)\Big(u_{\kappa,L}-\sum_{n=0}^N\phi^n_{j_1\ldots j_n}\nabla^n_{j_1\ldots j_n}\bar u_{\kappa,L,N}\Big)\\
\,=\,
\phi^N_{\kappa,j_1\ldots j_N}(\kappa^2-\partial_t^2)\nabla^N_{j_1\ldots j_N}\bar u_{\kappa,L,N}
+e_{j_{N+1}}\cdot A\big(\nabla\phi^N_{\kappa,j_1\ldots j_N}+\phi^{N-1}_{j_1\ldots j_{N-1}}e_{j_{N}}\big)\nabla^{N+1}_{j_1\ldots j_{N+1}}\bar u_{\kappa,L,N}\\
+\nabla\cdot (A\phi^N_{\kappa,j_1\ldots j_N}e_{j_{N+1}})\nabla^{N+1}_{j_1\ldots j_{N+1}}\bar u_{\kappa,L,N}
+e_{j_{N+2}}\cdot A\phi^N_{\kappa,j_1\ldots j_N}e_{j_{N+1}}\nabla^{N+2}_{j_1\ldots j_{N+2}}\bar u_{\kappa,L,N}\\
-\sum_{n=N+1}^{2N}\sum_{m=n-N}^{N}e_{j_{n}}\cdot\phi^{n-m-1}_{j_1\ldots j_{n-m-1}}\bar A^m_{j_{n-m}\ldots j_{n-2}}e_{j_{n-1}}\nabla^{n}_{j_1\ldots j_{n}}\bar u_{\kappa,L,N},
\end{multline*}
or equivalently, after recombining the terms and using the definition~\eqref{e.new-flux} 
of the flux correctors $\sigma^N_\kappa$ and $g^N_\kappa$,
\begin{multline}\label{eq:eqn-error}
\big(\partial_t^2-\nabla\cdot A\nabla\big)\Big(u_{\kappa,L}-\sum_{n=0}^N\phi^n_{j_1\ldots j_n}\nabla^n_{j_1\ldots j_n}\bar u_{\kappa,L,N}\Big)\\
\,=\,
\nabla\cdot \Big(\big(A\phi^N_{\kappa,j_1\ldots j_N}+(\sigma^N_{\kappa,j_1\ldots j_N})^T\big)\nabla\nabla^{N}_{j_1\ldots j_{N}}\bar u_{\kappa,L,N}\Big)
-\partial_t\Big(\phi^N_{\kappa,j_1\ldots j_N}\nabla^N_{j_1\ldots j_N}\partial_t\bar u_{\kappa,L,N}\Big)
\\
+\Big(\kappa^2 \phi^N_{\kappa,j_1\ldots j_n}+\kappa g^N_{\kappa,j_1\ldots j_N} \cdot \nabla-\sigma^N_{\kappa,j_1\ldots j_N}:\nabla^2\Big)\nabla^{N}_{j_1\ldots j_{N}}\bar u_{\kappa,L,N}\\
-\sum_{n=N+2}^{2N}\sum_{m=n-N}^{N}e_{j_{n}}\cdot\phi^{n-m-1}_{j_1\ldots j_{n-m-1}}\bar A^m_{j_{n-m}\ldots j_{n-2}}e_{j_{n-1}}\nabla^{n}_{j_1\ldots j_{n}}\bar u_{\kappa,L,N}.
\end{multline}
}
We now use the following standard estimate for the wave equation, see e.g.~\cite[Lemma~B.1]{DGR}: if $v$ satisfies
\[\left\{\begin{array}{l}
(\partial_t^2-\nabla\cdot A\nabla)v=\nabla\cdot f+\partial_tg+h,\\
v|_{t=0}=v^\circ,\\
\partial_tv|_{t=0}=0,
\end{array}\right.\]
 for some $f,g,h,v^\circ$ with $g|_{t=0}=0$, then we have
\begin{equation}\label{eq:apriori-wave-L2}
\|v^t\|_{\Ld^2(\R^d)}\,\lesssim\,\|v^\circ\|_{\Ld^2(\R^d)}+t\sup_{0\le s\le t}\|(f^s,g^s)\|_{\Ld^2(\R^d)}+t\sup_{0\le s\le t}\Big\|\int_0^sh\Big\|_{\Ld^2(\R^d)}.
\end{equation}
Applying this to the above wave equation~\eqref{eq:eqn-error},
we get
\begin{equation*}
\Big\|u_{\kappa,L}^t-\sum_{n=0}^N\phi^n_{j_1\ldots j_n}\nabla^n_{j_1\ldots j_n}\bar u_{\kappa,L,N}^t\Big\|_{\Ld^2(\R^d)}
\,\lesssim\,A_{\kappa,L,N}^\circ+t\big(B_{\kappa,L,N}^t+C_{\kappa,L,N}^t+D_{\kappa,L,N}^t\big),
\end{equation*}
and thus by the triangle inequality,
\begin{equation}\label{eq:estimL2-u}
\|u_{\kappa,L}^t-\bar u_{\kappa,L,N}^t\|_{\Ld^2(\R^d)}
\,\lesssim\,A_{\kappa,L,N}^\circ+A_{\kappa,L,N}^t+t\big(B_{\kappa,L,N}^t+C_{\kappa,L,N}^t+D_{\kappa,L,N}^t\big),
\end{equation}
where we have set for abbreviation
\begin{eqnarray}
A_{\kappa,L,N}^\circ&:=&\Big\|\sum_{n=1}^N\phi^n_{j_1\ldots j_n}\nabla^n_{j_1\ldots j_n}u^\circ_{\kappa,L}\Big\|_{\Ld^2(\R^d)},\label{eq:def-ABCD}\\
A_{\kappa,L,N}^t&:=&\Big\|\sum_{n=1}^N\phi^n_{j_1\ldots j_n}\nabla^n_{j_1\ldots j_n}\bar u_{\kappa,L,N}^t\Big\|_{\Ld^2(\R^d)},\nonumber\\
B_{\kappa,L,N}^t&:=&\sup_{0\le s\le t}\Big\|(\phi^N_\kappa,\sigma^N_\kappa)(\partial_s,\nabla)\nabla^{N}\bar u_{\kappa,L,N}^s\Big\|_{\Ld^2(\R^d)}\nonumber\\
C_{\kappa,L,N}^t&:=&\sup_{0\le s\le t}\Big\|\big(\kappa^2\phi^N_\kappa+\kappa g^N_\kappa\cdot \nabla-\sigma^N_\kappa :\nabla^2\big)\nabla^{N}\int_0^s\bar u_{\kappa,L,N}\Big\|_{\Ld^2(\R^d)},\nonumber\\
D_{\kappa,L,N}^t&:=&\sum_{n=N+2}^{2N}\sum_{m=n-N}^N\sup_{0\le s\le t}\Big\|\phi^{n-m-1}\bar A^m\nabla^{n}\int_0^s\bar u_{\kappa,L,N}\Big\|_{\Ld^2(\R^d)}.\nonumber
\end{eqnarray}

\medskip
\step4 Proof that, under assumption~\eqref{eq:cond-small-lambda}, we have for all $t\ge 0$,
\begin{eqnarray}
A_{\kappa,L,N}^t&\lesssim&(\kappa+K_N\kappa^{N-\delta})(L+\kappa^{-1}+t)^\e ,\label{eq:estim-A}\\
B_{\kappa,L,N}^t&\lesssim&K_N\kappa^{N+1-\delta}(L+\kappa^{-1}+t)^\e,\label{eq:estim-B}\\
C_{\kappa,L,N}^t&\lesssim&K_N\kappa^{N+1-\delta}(L+\kappa^{-1}+t)^\e ,\label{eq:estim-C}\\
D_{\kappa,L,N}^t&\lesssim&M_{N,\kappa}\kappa^{N+1}(L+\kappa^{-1}+t)^\e,\label{eq:estim-D}
\end{eqnarray}
where we have set for abbreviation
\begin{equation}\label{eq:abbr-LNl}
M_{N,\kappa}\,:=\,\sum_{n=1}^{N-1}\kappa^{n-1}\sum_{m=1}^{N-n}K_{N-m}K_{n+m},
\end{equation}
so that~\eqref{eq:estimL2-u} becomes
\begin{equation}\label{eq:estimL2-u-re}
\|u_{\kappa,L}^t-\bar u_{\kappa,L,N}^t\|_{\Ld^2(\R^d)}
\,\lesssim\,
\Big(\kappa+K_N\kappa^{N-\delta}
+(K_N+M_{N,\kappa})\kappa^{N+1-\delta}t\Big)(L+\kappa^{-1}+t)^{\e}.
\end{equation}
We start with the proof of~\eqref{eq:estim-A}. By definition~\eqref{eq:def-ABCD} of $A_{\kappa,L,N}$, the corrector bounds~\eqref{eq:estim-corr} allow to estimate
\begin{equation}\label{eq:estim-A-1}
A_{\kappa,L,N}^t\,\lesssim\,\sum_{n=1}^{N-1}K_n\|\langle\cdot\rangle^\e\nabla^n\bar u_{\kappa,L,N}\|_{\Ld^2(\R^d)}+K_N\kappa^{-\delta}\|\langle\cdot\rangle^\e\nabla^N\bar u_{\kappa,L,N}\|_{\Ld^2(\R^d)}.
\end{equation}
Using the Fourier formula~\eqref{eq:baru-lam} for $\bar u_{\kappa,L,N}$, together with~\eqref{eq:lower-symb}, we can easily compute
\begin{eqnarray*}
\|\nabla^n\bar u_{\kappa,L,N}^t\|_{\Ld^2(\R^d)}
&\le&\|\nabla^nu^\circ_{\kappa,L}\|_{\Ld^2(\R^d)},\\
\||\!\cdot\!|\nabla^n\bar u_{\kappa,L,N}^t\|_{\Ld^2(\R^d)}
&\lesssim&\|(|\!\cdot\!|+t)\nabla^nu^\circ_{\kappa,L}\|_{\Ld^2(\R^d)},
\end{eqnarray*}
and thus, by interpolation,
\begin{equation}\label{e.interp}
\|\langle\cdot\rangle^\e\nabla^n\bar u_{\kappa,L,N}^t\|_{\Ld^2(\R^d)}\,\lesssim\,\|(\langle\cdot\rangle+t)\nabla^nu^\circ_{\kappa,L}\|_{\Ld^2(\R^d)}^\e \|\nabla^nu^\circ_{\kappa,L}\|_{\Ld^2(\R^d)}^{1-\e}.
\end{equation}
Inserting this into~\eqref{eq:estim-A-1}, together with the definition~\eqref{eq:def-psilambda} of the regularized initial data~$u_{\kappa,L}^\circ$, we  get
\[A_{\kappa,L,N}^t\,\lesssim\,(L+\kappa^{-1}+t)^\e\sum_{n=1}^{N-1}K_n\kappa^{n}+(L+\kappa^{-1}+t)^\e K_N\kappa^{N-\delta},\]
and the claim~\eqref{eq:estim-A} follows under assumption~\eqref{eq:cond-small-lambda}.
A similar argument leads to~\eqref{eq:estim-B}.
We turn to the proof of~\eqref{eq:estim-C}. By definition~\eqref{eq:def-ABCD} of $C_{\kappa,L,N}$, the corrector bounds~\eqref{eq:estim-corr} allow to estimate
\[C_{\kappa,L,N}^t\,\lesssim\,K_N\kappa^{-\delta}\Big\|\langle\cdot\rangle^\e(\kappa^{2},\kappa\nabla,\nabla^2)\nabla^N\int_0^s\bar u_{\kappa,L,N}\Big\|_{\Ld^2(\R^d)}.\]
Noting that the Fourier formula~\eqref{eq:baru-lam} for $\bar u_{\kappa,L,N}$ yields
\[\int_0^s\nabla\bar u_{\kappa,L,N}\,=\,\int_{\R^d}e^{ix\cdot\xi}\tfrac{i\xi}{\sqrt{\mu_N(i\xi)}}\sin\big(s\sqrt{\mu_N(i\xi)}\big)\,\hat u^\circ_{\kappa,L}(\xi)\,d^*\xi,\]
a similar argument as above leads us to~\eqref{eq:estim-C}. Finally, for $D_{\kappa,L,N}$, a similar computation yields
\[D_{\kappa,L,N}^t\,\lesssim\,(L+\kappa^{-1}+t)^\e\sum_{n=N+2}^{2N}\kappa^{n-1}\sum_{m=n-N}^NK_{n-m-1}K_m,\]
and the claim~\eqref{eq:estim-D} follows in terms of the short-hand notation~\eqref{eq:abbr-LNl}.

\medskip
\step5 Dispersive estimates for the homogenized solution $\bar u_{\kappa,L,N}$.\\
We show that the solution $\bar u_{\kappa,L,N}$ of the complicated homogenized wave equation~\eqref{eq:0homog-baru} including dispersive corrections satisfies similar dispersive estimates as classical waves in an homogeneous medium, cf.~\eqref{eq:standard-disp-1} and~\eqref{eq:standard-disp-2}.
More precisely, we prove the following:
\begin{enumerate}[(i)]
\item Given $\eta>0$, provided that $\kappa$ is small enough in the sense that it satisfies the following strengthened version of~\eqref{eq:cond-small-lambda} for some $k\gg_\eta1$,
\begin{equation}\label{eq:cond-small-lambda/k=d}
\max_{0\le j\le k+1}\sum_{n=2\vee(j-1)}^N\tfrac{(n+1)!}{(n+1-j)!}K_n\kappa^{n-1}\,\le\,\tfrac1{2C_0},
\end{equation}
we have for all $R,L\ge1$ and $t\gg R+L$,
\begin{equation}\label{eq:estim-dispers}
\|\bar u_{\kappa,L,N}^t\|_{\Ld^\infty(B_R)}\,\lesssim\,\kappa^d(1+\kappa t)^{\eta-d}\|u^\circ\|_{\Ld^1(B_L)}.
\end{equation}
\item Provided that~\eqref{eq:cond-small-lambda/k=d} with $k=d-1$, we have for all $L\ge1$ and $t\ge0$,
\begin{equation}\label{eq:estim-dispers-2}
\|\bar u_{\kappa,L,N}^t\|_{\Ld^\infty(\R^d)}\,\lesssim\,\kappa^d(1+\kappa t)^{-\frac{d-1}2}\|u^\circ\|_{\Ld^1(\R^d)}.
\end{equation}
\end{enumerate}
As $\bar u_{\kappa,L,N}$ solves the homogenized wave equation~\eqref{eq:homog-baru} with initial data $u^\circ_{\kappa,L}=\rho_\kappa\ast(\chi_Lu^\circ)$, it can be represented as
\[\bar u_{\kappa,L,N}=G_{\kappa,N}\ast(\chi_Lu^\circ),\]
where $G_{\kappa,N}$ is the smoothened Green's function, solution of 
\begin{equation*}
\left\{\begin{array}{l}
\partial_t^2 G_{\kappa,N}=\nabla\cdot\big(\sum_{n=1}^N\bar A^n_{i_1\ldots i_{n-1}}\nabla^{n-1}_{i_1\ldots i_{n-1}}\big)\nabla G_{\kappa,N},\\
G_{\kappa,N}|_{t=0}=\rho_\kappa,\\
\partial_tG_{\kappa,N}|_{t=0}=0,
\end{array}\right.
\end{equation*}
which is understood as in~\eqref{eq:baru-lam} in terms of the Fourier formula
\begin{equation}\label{eq:smooth-green-four}
G_{\kappa,N}^t(x)\,:=\,\int_{\R^d}e^{ix\cdot\xi}\cos\big(t\sqrt{\mu_{N}(i\xi)}\big)\,\hat\rho_{\kappa}(\xi)\,d^*\xi.
\end{equation}
This representation of $\bar u_{\kappa,L,N}$ leads to the following local estimate, {for all $x\in B_R$,
\[|\bar u_{\kappa,L,N}(x)|
\,\le\,\int_{\R^d} |G_{\kappa,N}(y)||\chi_L(x-y) u^\circ(x-y)|\,dy\,\le\,\int_{B_{L+R}} |G_{\kappa,N}(y)||u^\circ(x-y)|\,dy,\]
and thus
\begin{equation}\label{eq:decomp-Linfty-GN}
\|\bar u_{\kappa,L,N}\|_{\Ld^\infty(B_R)}\,\le\,\|u^\circ\|_{\Ld^1(B_L)}\|G_{\kappa,N}\|_{\Ld^\infty(B_{R+L})}.
\end{equation}
Likewise,}
\begin{equation}
\|\bar u_{\kappa,L,N}\|_{\Ld^\infty(\R^d)}
\,\lesssim\,\|u^\circ\|_{\Ld^1(\R^d)}\|G_{\kappa,N}\|_{\Ld^\infty(\R^d)},\label{eq:decomp-Linfty-GN-re}
\end{equation}
so it remains to estimate the supremum norm of $G_{\kappa,N}$.
We start by using a $\kappa$-rescaling to rewrite~\eqref{eq:smooth-green-four} as
\begin{equation}\label{eq:smooth-green-four-re}
G_{\kappa,N}^t(x)\,=\,\kappa^d\tilde G_{\kappa,N}^{\kappa t}(\kappa x),\qquad\,\tilde G_{\kappa,N}^t(x):=\,\int_{\R^d}e^{i x\cdot\xi}\cos\big( t\mu_{\kappa,N}(i\xi)^{\frac12}\big)\,\hat\rho(\xi)\,d^*\xi,
\end{equation}
in terms of
\begin{equation}\label{eq:def-muNlam}
\mu_{\kappa,N}(i\xi)\,:=\,\kappa^{-2}\mu_{N}(i\kappa\xi)\,=\,\xi\cdot\Big(\sum_{n=1}^N\kappa^{n-1}(i\xi)^{\otimes(n-1)}_{j_1\ldots j_{n-1}}\bar A_{j_1\ldots j_{n-1}}^n\Big)\xi.
\end{equation}
Next, we decompose $\tilde G_{\kappa,N}$ as
\[\tilde G_{\kappa,N}\,=\,
\tfrac12\big(\tilde G_{\kappa,N,+}+\tilde G_{\kappa,N,-}\big),\]
where $\tilde G_{\kappa,N,+}$ and $\tilde G_{\kappa,N,-}$ are given by
\begin{equation*}
\tilde G_{\kappa,N,\pm}^t(x)\,:=\,\int_{\R^d}\exp\Big({i\big( x\cdot\xi\pm  t\mu_{\kappa,N}(i\xi)^\frac12\big)}\Big)\hat \rho(\xi)\,d^*\xi.
\end{equation*}
Given $0<\zeta<1$ to be later optimized, we further decompose
\begin{equation}\label{eq:GlamNpm}
\tilde G_{\kappa,N,\pm}^t(x)\,=\,\tilde G_{\kappa,N,\pm;\zeta,0}^t(x)+\tilde G_{\kappa,N,\pm;\zeta,1}^t(x),
\end{equation}
in terms of
\begin{gather}
\tilde G_{\kappa,N,\pm;\zeta,\sigma}^t(x)\,:=\,\int_{\R^d}\exp\Big({i\big( x\cdot\xi\pm  t\mu_{\kappa,N}(i\xi)^\frac12\big)}\Big)\hat \rho_{\zeta,\sigma}(\xi)\,d^*\xi,\label{eq:GlamNpm-def}\\
\hat \rho_{\zeta,0}(\xi)\,=\,\hat\rho(\xi)\chi(\zeta^{-1}\xi),
\qquad \hat \rho_{\zeta,1}(\xi)\,=\,\hat\rho(\xi)\big(1-\chi(\zeta^{-1}\xi)\big).\nonumber
\end{gather}
In order to estimate these oscillatory integrals, we argue as e.g.~in~\cite[Section~8.1.3]{BCD-11}, using that
\begin{multline*}
\exp\Big({i\big( x\cdot\xi\pm  t\mu_{\kappa,N}(i\xi)^\frac12\big)}\Big)\\
\,=\,\Big(1+t\big|\tfrac{x}{t}\pm\nu_{\kappa,N}(\xi)\big|^2\Big)^{-1}
\Big(1-i\big(\tfrac{x}{t}\pm\nu_{\kappa,N}(\xi)\big)\cdot\nabla_\xi\Big)\exp\Big({i\big( x\cdot\xi\pm  t\mu_{\kappa,N}(i\xi)^\frac12\big)}\Big),
\end{multline*}
with the short-hand notation
\begin{equation}\label{eq:def-nuN}
\nu_{\kappa,N}(\xi)\,:=\,\tfrac{\nabla_\xi\mu_{\kappa,N}(i\xi)}{2\sqrt{\mu_{\kappa,N}(i\xi)}}.
\end{equation}
Inserting this into~\eqref{eq:GlamNpm-def} and performing $k\ge0$ integrations by parts, we find
\begin{equation}\label{eq:IPP-D}
\tilde G_{\kappa,N,\pm;\zeta,\sigma}^t(x)\,=\,\int_{\R^d}\exp\Big({i\big( x\cdot\xi\pm  t\mu_{\kappa,N}(i\xi)^\frac12\big)}\Big)\big(D_{\xi;x,t,\pm}^{k}\hat \rho_{\zeta,\sigma}\big)(\xi)\,d^*\xi,
\end{equation}
provided that the integral makes sense,
where $D_{\xi;x,t,\pm}$ stands for the first-order differential operator
\begin{eqnarray}
\hspace{-1cm}D_{\xi;x,t,\pm}&:=&\Big(1+i\nabla_\xi\cdot\big(\tfrac xt\pm\nu_{\kappa,N}(\xi)\big)\Big)\Big(1+t\big|\tfrac{x}{t}\pm\nu_{\kappa,N}(\xi)\big|^2\Big)^{-1}\nonumber\\
&=&\tfrac{1\pm i(\nabla\cdot\nu_{\kappa,N})(\xi)}{1+t|\frac{x}{t}\pm\nu_{\kappa,N}(\xi)|^2}
\mp2t\tfrac{i(\frac xt\pm\nu_{\kappa,N}(\xi))^{\otimes2}:\nabla\nu_{\kappa,N}(\xi)}{(1+t|\frac{x}{t}\pm\nu_{\kappa,N}(\xi)|^2)^2}
+\tfrac{i(\frac xt\pm\nu_{\kappa,N}(\xi))}{1+t|\frac{x}{t}\pm\nu_{\kappa,N}(\xi)|^2}\cdot\nabla_\xi.\label{eq:def-Dxi}
\end{eqnarray}
 For all $k\ge0$, provided that $\kappa$ is small enough in the sense of~\eqref{eq:cond-small-lambda/k=d}, a similar argument as in~\eqref{eq:lower-symb} yields for all $\xi\in B\setminus\{0\}$ and $0\le j\le k$,
\begin{equation}\label{eq:estim-nuN}
|\nu_{\kappa,N}(\xi)|\,\simeq\,1,\qquad |\nabla^j_\xi\nu_{\kappa,N}(\xi)|\,\lesssim_j\,|\xi|^{-j}.
\end{equation}
A direct computation starting from~\eqref{eq:def-Dxi} then yields for all $\xi\ne0$ and $k\ge0$, provided that~\eqref{eq:cond-small-lambda/k=d} holds,
\begin{equation}\label{eq:Dkest0}
\big|\big(D_{\xi;x,t,\pm}^{k}\hat\rho_{\zeta,\sigma}\big)(\xi)\big|\,\lesssim_k\,\frac{\big(1+|\tfrac xt\pm\nu_{\kappa,N}(\xi)|\big)^{k-1}|\xi|^{-k}+|\tfrac xt\pm\nu_{\kappa,N}(\xi)|^k\zeta^{-k}}{(1+t|\tfrac xt\pm\nu_{\kappa,N}(\xi)|^2)^{k}}\,\mathds1_{|\xi|\le1}.
\end{equation}
We shall use this to bound the integral expression in~\eqref{eq:IPP-D}, and we split the proof into two further substeps, separately establishing~\eqref{eq:estim-dispers} and~\eqref{eq:estim-dispers-2}.

\medskip

\substep{5.1} Proof of~\eqref{eq:estim-dispers}.\\
For $|x|\ll t$, the bound~\eqref{eq:estim-nuN} ensures $|\frac xt\pm\nu_N(\xi)|^2\simeq1$ for $\xi\in B\setminus\{0\}$, so that~\eqref{eq:Dkest0} takes on the following guise: for all $\xi\ne0$ and $k\ge0$, provided that~\eqref{eq:cond-small-lambda/k=d} holds,
\begin{equation*}
\big|\big(D_{\xi;x,t,\pm}^{k}\hat\rho_{\zeta,\sigma}\big)(\xi)\big|\,\lesssim_k\,(1+t)^{-k}(|\xi|^{-k}+\zeta^{-k})\mathds1_{|\xi|\le1}.
\end{equation*}
Inserting this into~\eqref{eq:GlamNpm} and~\eqref{eq:IPP-D}, choosing $k=d-1$ in case $\sigma=0$, choosing any~$k\ge0$ in case $\sigma=1$, and using properties of $\rho,\chi$, we deduce
we deduce
\begin{eqnarray*}
\sup_{x:|x|\ll t}|\tilde G_{\kappa,N,\pm}^t(x)|&\lesssim&\int_{\R^d}\big|\big(D_{\xi;x,t,\pm}^{d-1}\hat \rho_{\zeta,0}\big)(\xi)\big|\,d\xi+\int_{\R^d}\big|\big(D_{\xi;x,t,\pm}^{k}\hat \rho_{\zeta,1}\big)(\xi)\big|\,d\xi\\
&\lesssim_k&(1+t)^{1-d}\int_{|\xi|\le\zeta}|\xi|^{1-d}d\xi+(1+t)^{-k}\zeta^{-k}\int_{\frac12\zeta\le|\xi|\le1}d\xi\\
&\lesssim&(1+t)^{1-d}\zeta+(1+t)^{-k}\zeta^{-k}.
\end{eqnarray*}
Optimizing with respect to $0<\zeta<\kappa$ and choosing $k>0$ sufficiently large, this yields the following result:
given $\eta>0$, provided that~\eqref{eq:cond-small-lambda/k=d} holds for some $k\gg_\eta1$ large enough, we have for all $t\ge0$,
\begin{equation}\label{eq:estim-G-disp1}
\sup_{x:|x|\ll t}|\tilde G_{\kappa,N,\pm}^t(x)|\,\lesssim_k\,(1+t)^{\eta-d},
\end{equation}
and therefore, by rescaling~\eqref{eq:smooth-green-four-re},
\begin{equation*}
\sup_{x:|x|\ll t}|G_{\kappa,N,\pm}^{t}(x)|\,\lesssim_k\,\kappa^d(1+\kappa t)^{\eta-d}.
\end{equation*}
Combined with~\eqref{eq:decomp-Linfty-GN}, this precisely yields the claim~\eqref{eq:estim-dispers}.

\medskip

\substep{5.2} Proof of~\eqref{eq:estim-dispers-2}.\\
We directly find from~\eqref{eq:smooth-green-four} that $\|G_{\kappa,N}^t\|_{\Ld^\infty(\R^d)}\lesssim\kappa^d$. After combination with~\eqref{eq:decomp-Linfty-GN-re}, this already yields the claim~\eqref{eq:estim-dispers-2} in dimension $d=1$, as well as in dimension $d>1$ for $\kappa t\le1$.
It remains to consider the higher-dimensional case $d>1$ for $\kappa t>1$.
Inserting~\eqref{eq:Dkest0} into~\eqref{eq:GlamNpm} and~\eqref{eq:IPP-D}, choosing $k=d-1$ both in cases~$\sigma=0$ and~$\sigma=1$, setting $\zeta=1$, and using properties of $\rho,\chi$,
we find for all~$x\in\R^d$ and~$t\ge0$,
\begin{multline}
|\tilde G_{\kappa,N,\pm}^t(x)|\,\lesssim\,\int_{\R^d}\big|\big(D_{\xi;x,t,\pm}^{d-1}\hat \rho_{\kappa;1,0}\big)(\xi)\big|\,d\xi+\int_{\R^d}\big|\big(D_{\xi;x,t,\pm}^{d-1}\hat \rho_{\kappa;1,1}\big)(\xi)\big|\,d\xi\nonumber\\
\,\lesssim_k\,\int_{|\xi|\le1}\frac{(1+|\tfrac xt\pm\nu_{\kappa,N}(\xi)|)^{d-2}|\xi|^{1-d}+|\tfrac xt\pm\nu_{\kappa,N}(\xi)|^{d-1}}{(1+t|\tfrac xt\pm\nu_{\kappa,N}(\xi)|^2)^{d-1}}\,d\xi.
\label{eq:estimGlNpm-diff}
\end{multline}
Noting that for all $t,A\ge0$,
\[\tfrac{A}{1+tA^2}\,\lesssim\,t^{-\frac12}(1+tA^2)^{-\frac12},\]
the above can be rewritten as
\begin{multline}\label{eq:estimGlNpm-diff00}
|\tilde G_{\kappa,N,\pm}^t(x)|\,\lesssim_k\,
\int_{|\xi|\le1}(1+t|\tfrac xt\pm\nu_{\kappa,N}(\xi)|^2)^{1-d}\,|\xi|^{1-d}\,d\xi\\
+t^{1-\frac{d}2}\int_{|\xi|\le1}(1+t|\tfrac xt\pm\nu_{\kappa,N}(\xi)|^2)^{-\frac{d}2}\,|\xi|^{1-d}\,d\xi\\
+t^{-\frac{d-1}2}\int_{|\xi|\le1}(1+t|\tfrac xt\pm\nu_{\kappa,N}(\xi)|^2)^{-\frac{d-1}2}\,d\xi.
\end{multline}
By the definition~\eqref{eq:def-muNlam} of $\mu_{\kappa,N}$, provided that~\eqref{eq:cond-small-lambda/k=d} holds for some~$k\ge1$, we can invert the map $H_{\kappa,N}:=\nabla \mu_{\kappa,N}(i\cdot):B\to\R^d$, and we find for all~$y\in H_{\kappa,N}(B)$,
\begin{equation}\label{eq:HN-unif}
|\nabla H_{\kappa,N}^{-1}(y)|\simeq1,\qquad|H_{\kappa,N}^{-1}(y)|\simeq|y|.
\end{equation}
Recalling the definition~\eqref{eq:def-nuN} of $\nu_{\kappa,N}$, the above estimate~\eqref{eq:estimGlNpm-diff00} then becomes after changing variables $y=H_{\kappa,N}(\xi)$,
\begin{multline}\label{eq:estimGlNpm-diff}
|\tilde G_{\kappa,N,\pm}^t(x)|\,\lesssim_k\,
\int_{|H_{\kappa,N}^{-1}(y)|\le1}(1+tA_{x,t,\pm}(y)^2)^{1-d}\,|y|^{1-d}\,dy\\
+t^{1-\frac{d}2}\int_{|H_{\kappa,N}^{-1}(y)|\le1}(1+tA_{x,t,\pm}(y)^2)^{-\frac{d}2}\,|y|^{1-d}\,dy\\
+t^{-\frac{d-1}2}\int_{|H_{\kappa,N}^{-1}(y)|\le1}(1+tA_{x,t,\pm}(y)^2)^{-\frac{d-1}2}\,dy,
\end{multline}
where we use the short-hand notation
\[A_{x,t,\pm}(y)\,:=\,\big|\tfrac xt\pm\tfrac12\mu_{\kappa,N}(iH_{\kappa,N}^{-1}(y))^{-\frac12}y\big|.\]
In order to evaluate these integrals, given $x\in\R^d$, we decompose the variable $y$ into its components $(y_1,y')$ parallel and orthogonal to $x$, respectively: more precisely, we define
\begin{equation*}
y\,:=\,y_1+y',\qquad y_1\,:=\,(\tfrac{ x\otimes x}{|x|^2})y,
\end{equation*}
and we note that we then have the lower bound
\begin{eqnarray*}
A_{x,t,\pm}(y)
&\ge&\big|(\Id-\tfrac{x\otimes x}{|x|^2})\tfrac12 \mu_{\kappa,N}(iH_{\kappa,N}^{-1}(y))^{-\frac12}y\big|\\
&=&\tfrac12 \mu_{\kappa,N}(iH_{\kappa,N}^{-1}(y))^{-\frac12}|y'|.
\end{eqnarray*}
Further using~\eqref{eq:lower-symb} and~\eqref{eq:HN-unif}, this yields
\begin{equation*}
A_{x,t,\pm}(y)
\,\gtrsim\, \tfrac{|y'|}{|y|}.
\end{equation*}
Inserting this into~\eqref{eq:estimGlNpm-diff} leads us to the following: provided~\eqref{eq:cond-small-lambda/k=d} holds for~$k=d-1$, we have for all $x\in\R^d$ and $t\ge0$,
\begin{multline*}
|\tilde G_{\kappa,N,\pm}^t(x)|\,\lesssim_k\,
\int_{|y|\le C}\Big(1+t(\tfrac{|y'|}{|y|})^2\Big)^{1-d}\,|y|^{1-d}\,dy\\
+t^{1-\frac{d}2}\int_{|y|\le C}\Big(1+t(\tfrac{|y'|}{|y|})^2\Big)^{-\frac{d}2}\,|y|^{1-d}\,dy
+t^{-\frac{d-1}2}\int_{|y|\le C}\Big(1+t(\tfrac{|y'|}{|y|})^2\Big)^{-\frac{d-1}2}\,dy.
\end{multline*}
Decomposing the integration domains according to $|y'|\ge\frac12|y|$ or $|y'|\le\frac12|y|$, and noting that in the latter case we have $|y'|\le|y_1|$ and $\frac23|y_1|\le|y|\le2|y_1|$, we can estimate
\begin{multline*}
|\tilde G_{\kappa,N,\pm}^t(x)|\,\lesssim_k\,
t^{1-d}
+\int_{|y'|\le|y_1|\le C}\Big(1+t(\tfrac{|y'|}{|y_1|})^2\Big)^{1-d}\,|y_1|^{1-d}\,dy\\
+t^{1-\frac{d}2}\int_{|y'|\le|y_1|\le C}\Big(1+t(\tfrac{|y'|}{|y_1|})^2\Big)^{-\frac{d}2}\,|y_1|^{1-d}\,dy\\
+t^{-\frac{d-1}2}\int_{|y'|\le|y_1|\le C}\Big(1+t(\tfrac{|y'|}{|y_1|})^2\Big)^{-\frac{d-1}2}\,dy.
\end{multline*}
The integrals are now easily evaluated and we get for all $t\ge1$,
\begin{equation*}
\|\tilde G_{\kappa,N,\pm}^t\|_{\Ld^\infty(\R^d)}\,\lesssim_k\,t^{-\frac{d-1}2},
\end{equation*}
and therefore, by rescaling~\eqref{eq:smooth-green-four-re}, for $\kappa t\ge1$,
\begin{equation*}
\|G_{\kappa,N,\pm}^{t}\|_{\Ld^\infty(\R^d)}\,\lesssim_k\,\kappa^d(\kappa t)^{-\frac{d-1}2}.
\end{equation*}
Combined with~\eqref{eq:decomp-Linfty-GN}, this precisely yields the claim~\eqref{eq:estim-dispers-2} for $\kappa t\ge1$, while the corresponding result for $\kappa t\le1$ is already proven.

\medskip
\step6 Conclusion.\\
By the triangle inequality, we can decompose for all $R,L,N\ge1$, {$0<\kappa\le 1$}, and $t\ge0$,
\[\|u^t\|_{\Ld^2(B_R)}\,\le\,\|u^t-u_{\kappa,L}^t\|_{\Ld^2(B_R)}+\|u_{\kappa,L}^t-\bar u_{\kappa,L,N}^t\|_{\Ld^2(\R^d)}+\|\bar u_{\kappa,L,N}^t\|_{\Ld^2(B_R)}.\]
Inserting~\eqref{eq:err-init-reg} and~\eqref{eq:estimL2-u-re}, we deduce for all $R,N\ge1$, $0<\kappa\le 1$, $t\ge0$, $L\ge 4(R+ct)$, and $n\ge1$, provided that~$\kappa$ is small enough in the sense of assumption~\eqref{eq:cond-small-lambda},
\begin{multline}\label{eq:estim-combin-disp0}
\|u^t\|_{\Ld^2(B_R)}\,\le\, C\kappa^{-1}\|\nabla u^\circ\|_{\Ld^2(\R^d)}
+\|\bar u_{\kappa,L,N}^t\|_{\Ld^2(B_R)}
+C_n(\kappa L)^{-n}\\
+C\Big(\kappa+K_N\kappa^{N-\delta}
+(K_N+M_{N,\kappa})\kappa^{N+1-\delta}t\Big)(L+\kappa^{-1}+t)^{\e}.
\end{multline}
If $u^\circ$ is already supported in $B_{L/2}$, we may use~\eqref{eq:err-init-reg+} instead of~\eqref{eq:err-init-reg}, and upgrade the above as follows: for all $R,N\ge1$, {$0<\kappa\le 1$}, and $t\ge 0$, provided that~$\kappa$ is small enough in the sense of assumption~\eqref{eq:cond-small-lambda},
\begin{multline}\label{eq:estim-combin-disp0+}
\|u^t\|_{\Ld^2(B_R)}\,\le\, C\kappa^{-1}\|\nabla u^\circ\|_{\Ld^2(\R^d)}
+\|\bar u_{\kappa,L,N}^t\|_{\Ld^2(B_R)}\\
+C\Big(\kappa+K_N\kappa^{N-\delta}+(K_N+M_{N,\kappa})\kappa^{N+1-\delta}t\Big)(L+\kappa^{-1}+t)^{\e}.
 \end{multline}
We split the rest of the proof into three further substeps, where we consider separately the periodic, quasiperiodic, and random settings.

\medskip
\substep{6.1} Periodic setting: proof of~(i).\\
In the periodic setting, by Lemma~\ref{lem:cor-per}, we can choose $N$ arbitrarily large, $\e=\delta=0$, and $K_n=C^n$.
Provided that $\kappa\ll1$ is small enough, letting $N\uparrow\infty$, the above result~\eqref{eq:estim-combin-disp0} then reads as follows: for all $R\ge1$, $t\ge0$,
$L\ge 4(R+ct)$, and $n\ge1$,
\begin{equation}\label{eq:per-predisp}
\|u^t\|_{\Ld^2(B_R)}\,\lesssim\,
\kappa^{-1}\|\nabla u^\circ\|_{\Ld^2(\R^d)}
+C_n(\kappa L)^{-n}
+\kappa
+\textstyle\liminf_{N\uparrow\infty}\|\bar u_{\kappa,L,N}^t\|_{\Ld^2(B_R)}.
\end{equation}
We now appeal to the dispersive estimates of Step~5 to bound the last right-hand side term.
By~\eqref{eq:estim-dispers-2}, provided that $\kappa\ll1$ is small enough, we get for all $R\ge1$, $t\ge0$, $L\ge 4(R+ct)$, and $n\ge1$,
\begin{equation*}
\|u^t\|_{\Ld^2(B_R)}\,\lesssim\,
\kappa^{-1}\|\nabla u^\circ\|_{\Ld^2(\R^d)}
+C_n(\kappa L)^{-n}+\kappa
+R^\frac d2\kappa^d(1+\kappa t)^{-\frac{d-1}2}\|u^\circ\|_{\Ld^1(\R^d)}.
\end{equation*}
Letting $L\uparrow\infty$, the second right-hand side term tends to $0$,
and the {first claimed dispersive estimate follows since   $\|u^\circ\|_{\Ld^2(\R^d)}=1$.}

\medskip\noindent
Next, if $u^\circ$ is supported in $B_{L/2}$, then by~\eqref{eq:estim-combin-disp0+} we may upgrade~\eqref{eq:per-predisp} as follows: for all~$R\ge1$ and $t\ge 0$, provided that $\kappa\ll1$ is small enough,
\begin{equation*} 
\|u^t\|_{\Ld^2(B_R)}\,\lesssim\,
\kappa^{-1}\|\nabla u^\circ\|_{\Ld^2(\R^d)}
+\kappa
+\textstyle\liminf_{N\uparrow\infty}\|\bar u_{\kappa,L,N}^t\|_{\Ld^2(B_R)}.
\end{equation*}
Now inserting the dispersive estimate~\eqref{eq:estim-dispers}, we get the following: given~$\eta>0$, provided that $\kappa\ll_\eta1$ is small enough, 
we have for all $R\ge1$ and~$t\gg R+L$,
\begin{equation*}
\|u^t\|_{\Ld^2(B_R)}\,\lesssim_\eta\,\kappa^{-1}\|\nabla u^\circ\|_{\Ld^2(\R^d)}
+\kappa
+R^\frac d2 \kappa^d(1+\kappa t)^{\eta-d}\|u^\circ\|_{\Ld^1(B_L)}.
\end{equation*}
To get the desired estimate for non-necessarily compactly supported initial conditions, it suffices to apply this result to $\chi_{L}u^\circ$ and to notice that
\begin{eqnarray*}
\|\nabla (\chi_{L}u^\circ)\|_{\Ld^2(\R^d)} &\lesssim& \|\nabla u^\circ\|_{\Ld^2(B_{L})}+L^{-1}\| u^\circ\|_{\Ld^2(B_{L})},\\
\|\chi_{L}u^\circ\|_{\Ld^p(\R^d)}&\lesssim& \|u^\circ\|_{\Ld^p(B_{L})}.
\end{eqnarray*}

\medskip
\substep{6.2} Quasiperiodic setting: proof of~(ii).\\
In the quasiperiodic setting, by Lemma~\ref{lem:cor-qper}, under a Diophantine condition and Gevrey regularity,
we can choose~$N$ arbitrarily large, $\e=\delta=0$, and $K_n=(Cn^{\theta})^n$ for some $\theta>0$.
The above result~\eqref{eq:estim-combin-disp0} then reads as follows:
for all $R,N\ge1$, $\kappa>0$, $t\ge0$, $L\ge4(R+ct)$, and $n\ge1$, provided that $\kappa N^\theta\ll1$, we have
\begin{multline}\label{eq:qper-predisp}
\|u^t\|_{\Ld^2(B_R)}\,\le\, C\kappa^{-1}\|\nabla u^\circ\|_{\Ld^2(\R^d)}
+\|\bar u_{\kappa,L,N}^t\|_{\Ld^2(B_R)}\\
+C_n(\kappa L)^{-n}
+C \kappa+(C\kappa N^\theta)^{ N}(1+\kappa t).
\end{multline}
We now appeal to the dispersive estimates of Step~5 to bound the second right-hand side term.
By~\eqref{eq:estim-dispers-2}, we get for all $R,N\ge1$, $\kappa>0$, $t\ge0$, $L\ge4(R+ct)$, and $n\ge1$, provided that~$\kappa N^\theta\ll1$,
\begin{multline*} 
\|u^t\|_{\Ld^2(B_R)}\,\le\, C\kappa^{-1}\|\nabla u^\circ\|_{\Ld^2(\R^d)}
+R^\frac d2\kappa^d(1+\kappa t)^{-\frac{d-1}2}\|u^\circ\|_{\Ld^1(\R^d)}\\
+C_n(\kappa L)^{-n}
+C \kappa+(C\kappa N^\theta)^{ N}(1+\kappa t).
\end{multline*}
Optimizing with respect to $N$ with $\kappa N^\theta\ll1$,
this becomes
\begin{multline*} 
\|u^t\|_{\Ld^2(B_R)}\,\le\, C\kappa^{-1}\|\nabla u^\circ\|_{\Ld^2(\R^d)}
+R^\frac d2\kappa^d(1+\kappa t)^{-\frac{d-1}2}\|u^\circ\|_{\Ld^1(\R^d)}\\
+C_n(\kappa L)^{-n}
+C\kappa\Big(1+t\exp(-\tfrac1C\kappa^{-\frac1\theta})\Big).
\end{multline*}
This yields the first part of item~(ii) by taking $t\le \exp(\tfrac1{C}\kappa^{-\frac1\theta})$ and letting $L\uparrow\infty$.

\medskip\noindent
Next, 
if $u^\circ$ is supported in $B_{L/2}$, then by~\eqref{eq:estim-combin-disp0+} we may upgrade~\eqref{eq:qper-predisp} as follows: for all~$R,L,N\ge1$, $\kappa>0$, and $t\ge 0$, provided that $\kappa N^\theta\ll1$ is small enough,
\[\|u^t\|_{\Ld^2(B_R)}\,\le\, C\kappa^{-1}\|\nabla u^\circ\|_{\Ld^2(\R^d)}+\|\bar u^t_{\kappa,L,N}\|_{\Ld^2(B_R)}+C\kappa+(C\kappa N^\theta)^N(1+\kappa t).\]
Now inserting the dispersive estimate~\eqref{eq:estim-dispers}, we get the following: given $\eta>0$, we obtain for all $R,L,N\ge1$, $\kappa>0$, and $t\gg R+L$, provided that $\kappa N^\theta\ll_\eta1$,
\begin{equation*}
\|u^t\|_{\Ld^2(B_R)}\,\lesssim_\eta\, \kappa^{-1}\|\nabla u^\circ\|_{\Ld^2(\R^d)}
+R^\frac d2\kappa^d(1+\kappa t)^{\eta-d}\|u^\circ\|_{\Ld^1(B_L)}
+\kappa
+(C\kappa N^\theta)^{ N}(1+\kappa t).
\end{equation*}
Optimizing with respect to $N$ as above, this entails
\begin{equation*}
\|u^t\|_{\Ld^2(B_R)}\,\lesssim_\eta\,\kappa^{-1}\|\nabla u^\circ\|_{\Ld^2(\R^d)}
+R^\frac d2\kappa^d(1+\kappa t)^{\eta-d}\|u^\circ\|_{\Ld^1(B_L)}\\
+\kappa\Big(1
+ t\exp(-\tfrac1C\kappa^{-\frac1\theta})\Big).
\end{equation*}
Applying this result to $\chi_Lu^\circ$ for non-necessarily compactly supported initial conditions~$u^\circ$,
this yields the second part of item~(ii).

\medskip\noindent
\substep{6.3} Random setting.\\
{In the random setting, by Lemma~\ref{lem:cor-rand} and the subsequent discussion in Section~\ref{sec:corr}, under suitable mixing assumptions,
we can choose $N=\lceil\frac d2\rceil$, any $\e>0$, $\delta=\frac12$ if~$d$ is odd, $\delta=0$ if~$d$ is even,
and $K_n\le \Kc_{\e}$ for some random variable $\Kc_{\e}$ with finite stretched exponential moments. 
Note that this yields $N-\delta=\frac d2$ whether $d$ is odd or even. 

\medskip\noindent
Set $\kappa_{*,\e}:=(\frac1C \Kc_{\e}^{-1})\wedge1$ with a constant $C$ large enough so that~\eqref{eq:cond-small-lambda} and~\eqref{eq:cond-small-lambda/k=d} with $k=d-1$ hold for~$\kappa\le\kappa_{*,\e}$.
The above result~\eqref{eq:estim-combin-disp0} then reads as follows: for all $R\ge1$, $0<\kappa\le\kappa_{*,\e}$, $t\ge0$, $L\ge4(R+ct)$, and $n\ge1$,
\begin{multline*}
\|u^t\|_{\Ld^2(B_R)}\,\lesssim\,\kappa^{-1}\|\nabla u^\circ\|_{\Ld^2(\R^d)}
+\|\bar u_{\kappa,L,N}^t\|_{\Ld^2(B_R)}\\
+C_n(\kappa L)^{-n}+\Big(\kappa+\kappa_{*,\e}^{-1}\kappa^{\frac d2}+\kappa_{*,\e}^{-2}\kappa^{\frac d2+1}t\Big)(L+\kappa^{-1}+t)^\e,
\end{multline*}
and thus, using the dispersive estimate~\eqref{eq:estim-dispers-2} of Step~5 to bound the second right-hand side term,
\begin{multline}
\|u^t\|_{\Ld^2(B_R)}\,\lesssim\, \kappa^{-1}\|\nabla u^\circ\|_{\Ld^2(\R^d)}
+R^{\frac d2}\kappa^d(1+\kappa t)^{-\frac{d-1}2}\|u^\circ\|_{\Ld^1(\R^d)}\\
+C_n(\kappa L)^{-n}+\Big(\kappa+\kappa_{*,\e}^{-1}\kappa^{\frac d2}+\kappa_{*,\e}^{-2}\kappa^{\frac d2+1}t\Big)(L+\kappa^{-1}+t)^\e,
\end{multline}
Restricting to times $t\le\kappa^{r-1-\frac d2}$ for some $0<r\le1$, and choosing $L\simeq\kappa^{r-1-\frac d2}$, we conclude for all $0<\kappa\le\kappa_{*,\e}$, $t\ge0$, $R\ge1$, provided that $R,t\le\kappa^{r-1-\frac d2}$,
\begin{equation*}
\|u^t\|_{\Ld^2(B_R)}\,\lesssim\,\kappa^{-1}\|\nabla u^\circ\|_{\Ld^2(\R^d)}
+R^{\frac d2}\kappa^d(1+\kappa t)^{-\frac{d-1}2}\|u^\circ\|_{\Ld^1(\R^d)}
+\kappa_{*,\e}^{-2}\kappa^{r-C\e}.
\end{equation*}
}
After redefining $\e$, and translating the origin,
this yields the first part of item~(iii).
The proof of the second part is similar and is omitted for shortness.
\qed

\section{Proof of the spreading of eigenstates in the lower spectrum}\label{sec:th}
This section is devoted to the proof of Theorem~\ref{th:main} and of Remark~\ref{rem:supp}.

\subsection{Proof of Theorem~\ref{th:main}}
{Assume that the heterogeneous acoustic operator~{$H_a=-\nabla\cdot A\nabla$} admits an eigenvalue~$\lambda>0$ with normalized eigenstate $\psi_\lambda\in H^1(\R^d)$, that is,
\begin{equation}\label{e.eigen}
-\nabla\cdot A\nabla\psi_\lambda=\lambda\psi_\lambda,\qquad\|\psi_\lambda\|_{\Ld^2(\R^d)}=1.
\end{equation}
Given {$0<\theta<\frac12$}, recall our definition~\eqref{eq:width} of the localization length of $\psi_\lambda$, that is, the minimal radius
\begin{equation}\label{eq:def-rho}
\ell_\theta(\psi_\lambda)\,:=\,\inf\big\{r\ge0:\|\psi_\lambda\|_{\Ld^2(B_r)}\ge1-\theta\big\} \,<\,\infty.
\end{equation}
To prove Theorem~\ref{th:main}, we take a dynamical approach and consider the large-scale behavior of the standing wave $(t,x)\mapsto\cos(t\sqrt\lambda)\psi_\lambda(x)$ associated with the eigenstate $\psi_\lambda$. As this eigenstate solves the heterogeneous wave equation, we can apply the large-scale dispersive estimates of Theorem~\ref{th:disp}.

\medskip

\step1 Periodic setting.
\\
Let $0<\eta<1$ be fixed. Write $U(t)\psi_\lambda=U(t)(\chi_{2L}\psi_\lambda)+U(t)((1- {\chi_{2L}})\psi_\lambda)$.
On the one hand, by unitarity,
\[\|U(t)((1- {\chi_{2L}})\psi_\lambda)\|_{\Ld^2(\R^d)} \le \|\psi_\lambda \|_{\Ld^2(\R^d \setminus B_L)}.\]
On the other hand, by Theorem~\ref{th:disp}(i), we have for all $0<\kappa \le\kappa_*$, $R,L\ge1$,  and $t\gg R+L$,
\begin{multline*}
\|U(t)( {\chi_{2L}}\psi_\lambda)\|_{\Ld^2(B_R)}\,
\lesssim \,
R^\frac d2 \kappa^d(1+\kappa t)^{\eta-d}\|\psi_\lambda\|_{\Ld^1({B_{2L}})}
\\+ \big(\kappa+(\kappa L)^{-1}\big)\|\psi_\lambda\|_{\Ld^2({B_{2L}})}+\kappa^{-1}\|\nabla \psi_\lambda\|_{\Ld^2({B_{2L}})}.
\end{multline*}
Hence, since $U(t)\psi_\lambda=\cos(t\sqrt\lambda) \psi_\lambda$, this entails for
$t \gg R+L$ with $t\in2\pi\sqrt\lambda{}^{-1}\N$ (which ensures $\cos(t\sqrt\lambda)=1$),
\begin{eqnarray*}
\|\psi_\lambda\|_{\Ld^2(B_R)} &\le& \|U(t)(  {\chi_{2L}}\psi_\lambda)\|_{\Ld^2(B_R)}
+\|U(t)( (1-{ \chi_{2L}})\psi_\lambda)\|_{\Ld^2(\R^d)}
\\
&\le & C R^\frac d2 \kappa^d(1+\kappa t)^{\eta-d}\|\psi_\lambda\|_{\Ld^1({B_{2L}})}
\\
&&+ C \big(\kappa+(\kappa L)^{-1}\big)\|\psi_\lambda\|_{\Ld^2({B_{2L}})}
+C\kappa^{-1}\|\nabla \psi_\lambda\|_{\Ld^2({B_{2L}})}
+\|\psi_\lambda \|_{\Ld^2(\R^d \setminus B_L)}.
\end{eqnarray*}
Letting $t\uparrow\infty$ to get rid of the first right-hand side term, and
then $L\uparrow\infty$, using the property $\int_{\R^d} \nabla \psi_\lambda \cdot A \nabla \psi_\lambda = \lambda$ of the eigenfunction, and choosing $R=\ell_\theta(\psi_\lambda)$,
this entails for all $0<\kappa\le\kappa_*$,
\begin{equation*}
0\,<\,{1-\theta}\,\le\,C\kappa+C\kappa^{-1}\sqrt\lambda .
\end{equation*}
Taking $\kappa:=\lambda^\frac14\le\kappa_*$, this yields a contradiction for $\lambda$ small enough. Hence $\ell_\theta(\psi_\lambda)=\infty$ and such an eigenfunction cannot exist.

\medskip
\step2 Quasiperiodic setting. 
\\
Let $0<\eta<1$ be fixed.
Proceeding as in the periodic setting, we obtain using Theorem~\ref{th:disp}(ii)
for all $0<\kappa \le\kappa_*$, $R,L\ge1$, and $R+L \ll t \le\exp(\tfrac1C\kappa^{-\frac1{\theta}})$ with  $t\in2\pi\sqrt\lambda{}^{-1}\N$,
\begin{eqnarray*}
\|\psi_\lambda\|_{\Ld^2(B_R)} &\le & C R^\frac d2 \kappa^d(1+\kappa t)^{\eta-d}\|\psi_\lambda\|_{\Ld^1({B_{2L}})}
\\
&&+ C\big(\kappa+(\kappa L)^{-1}\big)\|\psi_\lambda\|_{\Ld^2({B_{2L}})}
+C\kappa^{-1}\|\nabla \psi_\lambda\|_{\Ld^2({B_{2L}})}
+\|\psi_\lambda\|_{\Ld^2(\R^d \setminus {B_L})}.
\end{eqnarray*}
Choose $\kappa=\lambda^{\frac12-\e}\le\kappa_*$.
If $\ell_\theta(\psi_\lambda)\gtrsim \exp\big(\tfrac1C\lambda^{-(\frac12-\e)\frac1{\theta}}\big)$, there is nothing to prove. We may then assume that $\ell_\theta(\psi_\lambda)\ll \exp\big(\tfrac1C\lambda^{-(\frac12-\e)\frac1{\theta}}\big)$.
For the choice $R=\ell_\theta(\psi_\lambda)$, $L=\lambda^{-\frac12} \vee \ell_\theta(\psi_\lambda)$, and $t=\exp(\tfrac1C\lambda^{-(\frac12-\e)\frac1{\theta}})$, the above yields
\begin{multline*}
{1-2\theta}
\,\le\, C\ell_\theta(\psi_\lambda)^\frac d2 \lambda^{d(\frac12-\e)} (\lambda^{-\frac12} \vee \ell_\theta(\psi_\lambda) )^\frac d2  \lambda^{-(d-\eta)(\frac12-\e)}\exp \big(-\tfrac1C(d-\eta)\lambda^{-(\frac12-\e)\frac1{\theta}}\big)
\\
+ C(\lambda^{\frac12-\e}+\lambda^\e).
\end{multline*}
Absorbing the last right-hand side term into the left-hand side for $\lambda\ll_{\theta}1$ yields the claimed estimate.

\medskip
\step{3} Random setting.\\
Let $0<r,\eta,\e\le1$
be fixed.
Using now Theorem~\ref{th:disp}(iii), we obtain for all $0<\kappa \le\kappa_{*,\e}(0)$, $R,L\ge1$, and $R+L \ll t \le \kappa^{r-1-\frac d2}$,
\begin{multline*}
\|\psi_\lambda\|_{\Ld^2(B_R)} \,\le \, CR^\frac d2 \kappa^d(1+\kappa t)^{\eta-d}\|\psi_\lambda\|_{\Ld^1({B_{2L}})}
\\
+ C\big(\kappa_{*,\e}(0)^{-2}\kappa^{r-\e}+(\kappa L)^{-1}\big)\|\psi_\lambda\|_{\Ld^2({B_{2L}})}
+C\kappa^{-1}\|\nabla \psi_\lambda\|_{\Ld^2({B_{2L}})}
+\|\psi_\lambda\|_{\Ld^2(\R^d \setminus {B_L})}.
\end{multline*}
Choose $\kappa=\lambda^{\frac12-\e}\le\kappa_{*,\e}(0)$.
If $\ell_\theta(\psi_\lambda)\gtrsim \kappa^{r-1-\frac d2}=\lambda^{(\frac12-\e)(r- 1-\frac d2)}$, there is nothing to prove (upon redefining $\e$ in the statement). We may then assume that $\ell_\theta(\psi_\lambda)\ll \lambda^{(\frac12-\e)(r- 1-\frac d2)}$.
For the choice $R=\ell_\theta(\psi_\lambda)$, $L=\lambda^{-\frac12} \vee \ell_\theta(\psi_\lambda)$, and $t=\lambda^{(\frac12-\e) (r-1-\frac d2)}$, the above yields
\begin{multline*}
{1-2\theta}\,\le\, C\ell_\theta(\psi_\lambda)^\frac d2 \lambda^{d(\frac12-\e)} (\lambda^{-\frac12} \vee \ell_\theta(\psi_\lambda) )^\frac d2  \lambda^{-(d-\eta)(\frac12-\e)}
\lambda^{(d-\eta)(\frac12-\e)(\frac d2+1-r)}\\
+C\big(\kappa_{*,\e}(0)^{-2}\lambda^{(\frac12-\e)(r-\e)}+\lambda^\e\big).
\end{multline*}
 Absorbing the last right-hand side term into the left-hand side for $\lambda\ll_\theta \kappa_{*,\e}(0)^{\frac{4}{(1-2\e)(r-\e)}}$,
the conclusion follows by choosing $r$, $\e$, and $\eta$ small enough.}
\qed

\subsection{Proof of Remark~\ref{rem:supp}}\label{sec:remsupp}
The result is obtained as a post-processing of the result for the mass density using an exponentially-weighted energy estimate.
Given~$\alpha\ge1$, consider the exponential weight
\[\zeta_{\lambda,\alpha}(x)\,:=\,K_{\lambda,\alpha}^{-1}\exp(-\tfrac{\sqrt \lambda}{\alpha}|x|),\qquad K_{\lambda,\alpha}\,:=\,\int_{\R^d}\exp(-\tfrac{\sqrt \lambda}{\alpha}|x|)\,dx.\]
For any $z\in\R^d$, testing the eigenvalue relation~\eqref{e.eigen} with $\zeta_{\lambda,\alpha}(\cdot-z)\psi_\lambda$, we find
\begin{multline*}
\lambda \int_{\R^d} \zeta_{\lambda,\alpha}(\cdot-z)|\psi_\lambda|^2
\,=\,\int_{\R^d}\zeta_{\lambda,\alpha}(\cdot-z)\nabla \psi_\lambda\cdot A\nabla \psi_\lambda\\
-\tfrac{\sqrt\lambda}\alpha\int_{\R^d} \zeta_{\lambda,\alpha}(z-x) \psi_\lambda(x) \tfrac{z-x}{|z-x|} \cdot A \nabla \psi_\lambda(x)\,dx.
\end{multline*}
By Young's inequality, with $\alpha\ge1$, we deduce
\begin{eqnarray}
\lambda \int_{\R^d} \zeta_{\lambda,\alpha}(\cdot-z)|\psi_\lambda|^2
&\le&\tfrac{1+\frac1{2\alpha}}{1-\frac1{2\alpha}}\int_{\R^d}\zeta_{\lambda,\alpha}(\cdot-z)\nabla \psi_\lambda\cdot A\nabla \psi_\lambda\nonumber\\
&\le&(1+\tfrac2{\alpha})\int_{\R^d}\zeta_{\lambda,\alpha}(\cdot-z)\nabla \psi_\lambda\cdot A\nabla \psi_\lambda.\label{e.averaging}
\end{eqnarray}
Next, we upgrade this estimate to get a comparison between the local norms $\lambda \int_{B_R}|\psi_\lambda|^2$ and~\mbox{$\int_{B_R}\nabla \psi_\lambda\cdot A\nabla \psi_\lambda$} (with non-smooth cut-off).
For that purpose, using local averaging with the weight $\zeta_{\lambda,\alpha}$, we start by writing
\begin{equation}\label{eq:tosplit-average}
\lambda\int_{\R^d \setminus B_R}|\psi_\lambda |^2\,=\,\lambda\int_{\R^d}\Big(\int_{\R^d \setminus B_R}\zeta_{\lambda,\alpha}(x-z)|\psi_\lambda(x)|^2\,dx\Big)\,dz,
\end{equation}
and we separately consider the near-field and far-field contributions in this integral.
On the one hand, for the integral over $|z|\le R/2$, using that the conditions $|z|\le R/2$ and $|x|\ge R$ imply $|x-z|\ge R/2$, and using the normalization of $\psi_\lambda$, we find
\begin{eqnarray*}
\lambda\int_{B_{R/2}}\Big(\int_{\R^d \setminus B_R}\zeta_{\lambda,\alpha}(x-z)|\psi_\lambda(x)|^2\,dx\Big)\,dz
&\le&C\lambda R^d\|\zeta_{\lambda,\alpha}\|_{\Ld^\infty(\R^d\setminus B_{R/2})}\\
&\le&C\lambda\exp(-\tfrac{\sqrt\lambda}{4\alpha}R).
\end{eqnarray*}
On the other hand, for the integral over $|z|\ge R/2$, we appeal to~\eqref{e.averaging} and obtain
\begin{eqnarray*}
\lefteqn{\lambda\int_{\R^d \setminus B_{R/2}}\Big(\int_{\R^d \setminus B_R}\zeta_{\lambda,\alpha}(x-z)|\psi_\lambda(x)|^2\,dx\Big)\,dz}\\
&\le&\lambda\int_{\R^d \setminus B_{R/2}}\Big(\int_{\R^d}\zeta_{\lambda,\alpha}(\cdot-z)|\psi_\lambda|^2\Big)\,dz\\
&\le&(1+\tfrac2\alpha)\int_{\R^d \setminus B_{R/2}}\Big(\int_{\R^d}\zeta_{\lambda,\alpha}(\cdot-z)\nabla\psi_\lambda\cdot A\nabla\psi_\lambda\Big)\,dz.
\end{eqnarray*}
We further split the last integral in bracket over $\R^d\setminus B_{R/4}$ and $B_{R/4}$, and we note respectively that
\[\int_{\R^d \setminus B_{R/2}}\Big(\int_{\R^d\setminus B_{R/4}}\zeta_{\lambda,\alpha}(\cdot-z)\nabla\psi_\lambda\cdot A\nabla\psi_\lambda\Big)\,dz\,\le\,\int_{\R^d\setminus B_{R/4}}\nabla\psi_\lambda\cdot A\nabla\psi_\lambda,\]
and that
\begin{eqnarray*}
\lefteqn{\int_{\R^d \setminus B_{R/2}}\Big(\int_{B_{R/4}}\zeta_{\lambda,\alpha}(\cdot-z)\nabla\psi_\lambda\cdot A\nabla\psi_\lambda\Big)\,dz}\\
&\le&\Big(\int_{\R^d \setminus B_{R/2}}\|\zeta_{\lambda,\alpha}\|_{\Ld^\infty(B_{R/4}(z))}\,dz\Big)\Big(\int_{\R^d}\nabla\psi_\lambda\cdot A\nabla\psi_\lambda\Big)\\
&\le&C\lambda \exp(-\tfrac{\sqrt\lambda}{8\alpha}R).
\end{eqnarray*}
Combining these different estimates back into~\eqref{eq:tosplit-average}, we obtain
\begin{equation}\label{eq:tosplit-average-re}
\lambda\int_{\R^d \setminus B_{R}}|\psi_\lambda(x)|^2\,dx\,\le\,(1+\tfrac2\alpha)\int_{\R^d\setminus B_{R/4}}\nabla\psi_\lambda\cdot A\nabla\psi_\lambda+C\lambda \exp(-\tfrac{\sqrt\lambda}{8\alpha}R).
\end{equation}
We can now use this bound to compare the minimal radii $\ell_\theta(\psi_\lambda)$ and $\ell'_\theta(\psi_\lambda)$ defined in~\eqref{eq:width} and~\eqref{eq:width-bis}: for the choice $R:=\ell_\theta(\psi_\lambda)$, the left-hand side in~\eqref{eq:tosplit-average-re} is by definition equal to~$\lambda\theta^2$, so we can deduce
\[\tfrac1{\sqrt{\lambda}}\|\sqrt A\nabla\psi_\lambda\|_{\Ld^2(\R^d\setminus B_{R/4})}\,\ge\,\big(1+\tfrac2\alpha\big)^{-\frac12}\big(\theta^2-C\exp(-\tfrac{\sqrt\lambda}{8\alpha}R)\big)^\frac12,\]
which means
\[\ell'_{\hat \theta}(\psi_\lambda)\,\ge\,\tfrac14\ell_\theta(\psi_\lambda),\qquad\text{provided $\hat\theta\le\big(1+\tfrac2\alpha\big)^{-\frac12}\big(\theta^2-C\exp(-\tfrac{\sqrt\lambda}{8\alpha}\ell_\theta(\psi_\lambda))\big)^\frac12$}.\]
In all the cases considered, as stated in Theorem~\ref{th:main}, for any $\e>0$ and {$0<\theta<\frac12$}, we have shown that we have at least $\ell_\theta(\psi_\lambda)\ge\lambda^{\e-3/4}$ for $\lambda$ small enough, so the above becomes
\[\ell'_{\hat \theta}(\psi_\lambda)\,\ge\,\tfrac14\ell_\theta(\psi_\lambda),\qquad\text{provided $\hat\theta\le\big(1+\tfrac2\alpha\big)^{-\frac12}\big(\theta^2-C\exp(-\tfrac{\lambda^{-1/8}}{8\alpha})\big)^\frac12$}.\]
Choosing $\alpha$ arbitrarily large to reach any {$0<\hat\theta<\frac12$} up to taking $\lambda$ small enough, the statement of Remark~\ref{rem:supp} follows.
\qed

\endgroup
\section*{Acknowledgements}
{We thank Joshua McGinnis and Doug Wright for discussions on their recent work~\cite{McGinnis-Wright-23}, which gave the idea of using massive correctors in Section~\ref{sec:corr}.}
MD also thanks Fran\c{c}ois Pagano for preliminary discussions on the topic of this work,
and acknowledges financial support from the F.R.S.-FNRS, as well as from the European Union (ERC, PASTIS, Grant Agreement n$^\circ$101075879).
AG acknowledges financial support from the European Research Council (ERC) under the European Union's Horizon 2020 research and innovation programme (Grant Agreement n$^\circ$864066).\footnote{ Views and opinions expressed are however those of the authors only and do not necessarily reflect those of the European Union or the European Research Council Executive Agency. Neither the European Union nor the granting authority can be held responsible for them.}

\bibliographystyle{abbrv}
\bibliography{biblio}

\end{document}